\documentclass[12pt]{iopart}
\usepackage{amsthm,amssymb,amsfonts,amsbsy,latexsym}
\usepackage{graphicx,psfrag}
\usepackage{times}
\usepackage{color}

\newcommand{\Cset}{\mathbb{C}}
\newcommand{\Hset}{\mathbb{H}}
\newcommand{\Nset}{\mathbb{N}}

\newcommand{\Rset}{\mathbb{R}}
\newcommand{\Sset}{\mathbb{S}}
\newcommand{\Tset}{\mathbb{T}}
\newcommand{\Zset}{\mathbb{Z}}

\theoremstyle{plain}
\newtheorem{thm}{Theorem}
\newtheorem{pro}[thm]{Proposition}
\newtheorem{lem}[thm]{Lemma}

\newtheorem{defi}{Definition}
\theoremstyle{remark}
\newtheorem{remark}{Remark}

\newcommand{\rmD}{{\rm D}}
\newcommand{\Identity}{{\rm I}}
\newcommand{\Lipschitz}{\mathop{\rm Lip}\nolimits}
\newcommand{\Fnorm}[1]{ \| #1  \|}
\newcommand{\Snorm}[1]{ | #1  |}

\newcommand{\Length}{\mathop{\rm Length}\nolimits}
\newcommand{\Area}{\mathop{\rm Area}\nolimits}

\newcommand{\Graph}{\mathop{\rm graph}\nolimits}

\newif\iffigures
\figurestrue

\begin{document}

\title[On the length and area spectrum of analytic convex domains]
      {On the length and area spectrum of analytic convex domains}

\author{Pau Mart{\'\i}n\dag, Rafael Ram\'{\i}rez-Ros\ddag and Anna Tamarit-Sariol\ddag}

\address{\dag Departament de Matem\`{a}tica Aplicada IV,
         Universitat Polit\`{e}cnica de Catalunya,
         Ed.~C3, Jordi Girona 1--3, 08034 Barcelona, Spain}

\address{\ddag Departament de Matem\`{a}tica Aplicada I,
         Universitat Polit\`{e}cnica de Catalunya,
         Diagonal 647, 08028 Barcelona, Spain}

\eads{\mailto{Martin@ma4.upc.edu},
      \mailto{Rafael.Ramirez@upc.edu},
      \mailto{Anna.Tamarit@upc.edu}}

\begin{abstract}
Area-preserving twist maps have at least two different $(p,q)$-periodic orbits
and every $(p,q)$-periodic orbit has its $(p,q)$-periodic action
for suitable couples $(p,q)$.
We establish an exponentially small upper bound for the
differences of $(p,q)$-periodic actions
when the map is analytic on a $(m,n)$-resonant rotational invariant curve
(resonant RIC)
and $p/q$ is ``sufficiently close'' to $m/n$.
The exponent in this upper bound is closely related to the
analyticity strip width of a suitable angular variable.
The result is obtained in two steps.
First, we prove a Neishtadt-like theorem,
in which the $n$-th power of the twist map is written as
an integrable twist map plus an exponentially small remainder
on the distance to the RIC.
Second, we apply the MacKay-Meiss-Percival action principle.

We apply our exponentially small upper bound to several billiard problems.
The resonant RIC is a boundary of the phase space in almost all of them.
For instance, we show that the lengths (respectively, areas) of
all the $(1,q)$-periodic billiard (respectively, dual billiard) trajectories
inside (respectively, outside) analytic strictly convex domains
are exponentially close in the period $q$.
This improves some classical results of Marvizi, Melrose,
Colin de Verdi\`ere, Tabachnikov, and others about the smooth case.
\end{abstract}




\section{Introduction}

Billiards were introduced by Birkhoff~\cite{Birkhoff1966}.
Let $\Gamma$ be a smooth strictly convex curve in the plane,
oriented counterclockwise,
and let $\Omega$ be the billiard table enclosed by $\Gamma$.
Billiard trajectories inside $\Omega$ consist of polygonal lines
inscribed in $\Gamma$ whose consecutive sides obey to the rule
``the angle of reflection is equal to the angle of incidence.''
See~\cite{KatokHasselblatt1995,KozlovTreschev1991,Tabachnikov1995_Billiards}
for a general description.

A $(p,q)$-periodic billiard trajectory forms a closed
polygon with $q$ sides that makes $p$ turns inside $\Gamma$.
Birkhoff~\cite{Birkhoff1966} proved that there are at least two different
Birkhoff $(p,q)$-periodic billiard trajectories inside $\Omega$
for any relatively prime integers $p$ and $q$ such that $1 \le p \le q$.

Let $\mathcal{L}^{(p,q)}$ be the supremum of the absolute values of the
differences of the lengths of all such trajectories.
The quantities $\mathcal{L}^{(p,q)}$ were already studied by
Marvizi and Melrose~\cite{MarviziMelrose1982} and
Colin de Verdi\`ere~\cite{Colin1984} for smooth tables.
The former authors produced an asymptotic expansion of the lengths for
$(p,q)$-periodic billiard trajectories approaching $\Gamma$
when $p$ is fixed and $q \to +\infty$.
They saw that there exists a sequence $(l_k)_{k\ge1}$,
depending only on $p$ and $\Gamma$, such that,
if $L^{(p,q)}$ is the length of any $(p,q)$-periodic trajectory, then
\[
L^{(p,q)} \asymp
p \Length[\Gamma] + \sum_{k \ge 1} \frac{l_k}{q^{2k}},
\qquad q\to\infty,
\]
where $l_1 = l_1(\Gamma,p) = -
\frac{1}{24} \left(p \int_\Gamma \kappa^{2/3}(s) \rmd s \right)^3$,
and $\kappa(s)$ is the curvature of $\Gamma$ as a function of
the arc-length parameter $s$.
The symbol $\asymp$ means that the series in the right hand side
is asymptotic to $L^{(p_,q)}$.
The assymptotic coefficients $l_k=l_k(\Gamma,p)$ can be explicitly
written in terms of the curvature $\kappa(s)$.
For instance, the explicit formulas for $l_1$, $l_2$, $l_3$, and $l_4$ can
be found in~\cite{Sorrentino2014}.
Since the expansion of the lengths in powers of $q^{-1}$ coincides for
all these $(p,q)$-periodic trajectories,
$\mathcal{L}^{(p,q)}=\Or(q^{-\infty})$ for smooth strictly convex tables
when $p$ is fixed and $q \to +\infty$.
Colin de Verdi\`ere studied the lengths of periodic trajectories
close to an elliptic $(1,2)$-periodic trajectory on a smooth symmetric
billiard table,
and found that the quantities $\mathcal{L}^{(p,q)}$ are again
beyond all order with respect to $q$.

These works suggest that the supremum length differences $\mathcal{L}^{(p,q)}$ are
exponentially small in the period $q$ for analytic strictly convex tables.
Indeed, we have proved that if $\Gamma$ is analytic and $p$ is a fixed
positive integer, then there exists $K, q_*, \alpha > 0$ such that
\begin{equation}\label{eq:BoundDifferenceLengths}
\mathcal{L}^{(p,q)}\le K\rme^{-2\pi \alpha q/p},
\end{equation}
for all integer $q\ge q_*$ relatively prime with $p$.
The exponent $\alpha$ is related to the width of a complex strip where a certain
$1$-periodic angular coordinate is analytic.
A more precise statement is given in Theorem~\ref{thm:ClassicalBoundary}.

Similar exponentially small upper bounds hold in other billiard problems.
We mention two examples.
First, for $(p,q)$-periodic billiard trajectories inside
strictly convex analytic tables of constant width when $p/q\to 1/2$.
Second, for $(p,q)$-periodic billiard trajectories inside strictly convex
analytic tables in surfaces of constant curvature when $p/q\to 0$.

The billiard dynamics close to the boundary has also been studied
from the point of view of KAM theory.
Lazutkin~\cite{Lazutkin1973} proved that there are infinitely many
caustics inside any $C^{555}$ strictly convex table.
These caustics accumulate at the boundary of the table,
and have Diophantine rotation numbers.
Douady~\cite{Douady1982} improved the result to $C^7$ billiard tables.

A special remark on the relevance of these results is the following.
Kac~\cite{Kac1966} formulated the inverse spectral problem for planar domains.
That is, to study how much geometric information about $\Omega$
can be obtained from the Laplacian spectrum with homogeneous Dirichlet
conditions on $\Gamma$.
Andersson and Melrose~\cite{AnderssonMelrose1977} gave an explicit
relation between the length spectrum and the Laplacian spectrum.
The length spectrum of $\Omega$ is the union of the lengths of
all its $(p,q)$-periodic billiard trajectories and
all the integer multiples of $\Length[\Gamma]$.
See also~\cite{MarviziMelrose1982,Colin1984}.

Our results also apply to the dual billiards
introduced by Day~\cite{Day1947} and popularized by
Moser~\cite{Moser1978} as a crude model for planetary motion.
Some general references
are~\cite{GutkinKatok1995,Boyland1996,Tabachnikov1995,Tabachnikov1995_Billiards}.
Let $\mho$ be unbounded component of $\Rset^2 \setminus\Gamma$.
The dual billiard map $f: \mho \to \mho$ is defined as follows:
$f(z)$ is the reflection of $z$ in the tangency point of the
oriented tangent line to $\Gamma$ through $z$.
Billiards and dual billiards are projective dual in
the sphere~\cite{Tabachnikov1995}.

A $(p,q)$-periodic dual billiard trajectory forms a closed
circumscribed polygon with $q$ sides that makes $p$ turns outside $\Gamma$.
The area of a $(p,q)$-periodic trajectory is the area enclosed by the
corresponding polygon, taking into account some multiplicities if $p \ge 2$.
There are at least two different Birkhoff $(p,q)$-periodic dual billiard
trajectories outside $\Gamma$ for any relatively prime integers
$p$ and $q$ such that $q \ge 3$ and $1 \le p \le q$.

Tabachnikov~\cite{Tabachnikov1995_Billiards,Tabachnikov1995} studied
the supremum $\mathcal{A}^{(1,q)}$ of the absolute value of the
differences of the areas enclosed by all such $(1,q)$-periodic trajectories
for smooth tables.
He proved that there is a sequence $(a_k)_{k\ge 1}$,
depending only on $\Gamma$, such that,
if $A^{(1,q)}$ is the area enclosed by any $(1,q)$-periodic
dual billiard trajectory, then
\begin{equation}\label{eq:AsymptoticArea}
A^{(1,q)} \asymp
\Area[\Omega] + \sum_{k \ge 1} \frac{a_k}{q^{2k}},
\qquad q\to\infty,
\end{equation}
where $a_1 = a_1(\Gamma) = \frac{1}{24}\int_{\Gamma} \kappa^{1/3}(s) \rmd s$.
Hence, the expansion of the areas in powers of $q^{-1}$ coincides for
all these $(1,q)$-periodic trajectories, and so,
$\mathcal{A}^{(1,q)}=\Or(q^{-\infty})$ for smooth strictly
convex dual tables when $q \to +\infty$.
Douady~\cite{Douady1982} found the existence of infinitely many
invariant curves outside any  $C^7$ strictly convex dual table.
These invariant curves accumulate at the boundary of the dual table
and have Diophantine rotation numbers.

In a completely analogous way to (classical) billiards,
we have proved that, once fixed any positive integer $p$,
if $\Gamma$ is analytic, then there exists $K, q_*, \alpha > 0$ such that
\begin{equation}\label{eq:BoundDifferenceAreas}
\mathcal{A}^{(p,q)} \le K\rme^{-2\pi \alpha q/p},
\end{equation}
for all integer $q\ge q_*$ relatively prime with $p$.
Once more,
the exponent $\alpha$ is related to the width of a complex strip where a certain
$1$-periodic angular coordinate is analytic.
The precise statement is given in Theorem~\ref{thm:DualBoundary}.

Still in the context of dual billiards,
the points at infinity can be seen as $(1,2)$-periodic points,
hence they form a $(1,2)$-resonant RIC.
Douady~\cite{Douady1982} found the existence of infinitely many
invariant curves outside any $C^8$ strictly convex dual table.
These invariant curves accumulate at infinity and have Diophantine
rotation numbers.
We have proved that, once fixed any constant $L \ge 1$,
if $\Gamma$ is analytic, then there exist  $K, q_*, \alpha > 0$ such that
\begin{equation}\label{eq:BoundDifferenceAreas_bis}
\mathcal{A}^{(p,q)}\le
K \exp\left(-\frac{\pi \alpha}{|p/q-1/2|}\right),
\end{equation}
for all relatively prime integers $p$ and $q$ such that
$1 \le |2p-q|\le L$ and $q\ge q_*$.
See Theorem~\ref{thm:DualInfinity}.

The three exponents $\alpha$ that appear in the exponentially
small upper bounds~(\ref{eq:BoundDifferenceLengths}),
(\ref{eq:BoundDifferenceAreas}), and~(\ref{eq:BoundDifferenceAreas_bis})
may be different,
since each one is associated to a different analyticity strip width.
Besides, all of these upper bounds follow directly from a general
upper bound about analytic area-preserving twist maps with
analytic resonant RICs.
Let us explain it.

Classical and dual billiard maps are exact twist maps
defined on an open cylinder when written in suitable coordinates.
Exact twist maps have been vastly studied.
They satisfy a Lagrangian formulation and their orbits are stationary
points of the action functional.
See for instance~\cite{Birkhoff1966,Meiss1992,KatokHasselblatt1995}.

Birkhoff~\cite{Birkhoff1966} showed that the minima and minimax points of
the  $(p,q)$-periodic action correspond to two different
Birkhoff $(p,q)$-periodic orbits of the twist map.
A Birkhoff $(p,q)$-periodic orbit is an orbit such that, after $q$ iterates,
performs exactly $p$ revolutions around the cylinder and its points are
ordered in the base $\Tset$ as the ones following a rigid rotation of angle $p/q$.
Since there exist at least two different Birkhoff $(p,q)$-periodic orbits,
we consider the supremum $\Delta^{(p,q)}$ of the absolute value
of the differences of the actions among all of them.
The quantity $\Delta^{(p,q)}$ coincides with $\mathcal{L}^{(p,q)}$
and $\mathcal{A}^{(p,q)}$ for classical and dual billiards, respectively.

Let $\Delta W_{p/q}$ be the difference of actions between
the minimax and minima $(p,q)$-periodic orbits.
Note that $\Delta^{(p,q)}$ is an upper bound of  $\Delta W_{p/q}$.
Mather~\cite{Mather1986} used  $\Delta W_{p/q}$ as a criterion to prove the
existence of RICs of given irrational rotation numbers.
More concretely,
he proved that there exists a RIC with irrational rotation number $\varrho$
if and only if $\lim_{p/q \to \varrho} \Delta W_{p/q} = 0$.

Another criterion related to the destruction of RICs,
in this case empirical, was proposed by Greene.
The destruction of a RIC with Diophantine rotation number $\varrho$
under perturbation is related to a
``sudden change from stability to instability of the nearby periodic
orbits''~\cite{Greene1979}.
The stability of a periodic orbit is measured by the residue.
MacKay~\cite{MacKay1992} proved the criterion in some contexts.
In particular, for an analytic area-preserving twist map,
the residue of a sequence of periodic orbits with rotation numbers tending
to $\varrho$ decays exponentially in a positive power of the distance
between rotation numbers, that is, $|\varrho-p/q|^d$ for some~$d>0$.
The same proof leads to a similar exponentially small bound of
Mather's~$\Delta W_{p/q}$ as $p/q \to \varrho$.
Delshams and de la Llave~\cite{DelshamsLlave2000} studied similar
problems for analytic area-preserving nontwist maps.

Generically RICs with a rational rotation number break
under perturbation~\cite{RamirezRos2006,PintodeCarvalhoRamirezRos2013}.
Nevertheless, there are situations in which some distinguished
resonant RICs always exist.
See Sections~\ref{sec:LengthSpectrum} and~\ref{sec:AreaSpectrum} for several examples related
to billiard and dual billiard maps.

Let us assume that we have an analytic exact twist map
with a $(m,n)$-resonant RIC.
That is, a RIC whose points are $(m,n)$-periodic.
Then there exist some variables $(x,y)$ in which the resonant RIC
is located at $\{ y = 0\}$ and the $n$-th power of the exact twist map
is a small perturbation of the integrable twist map $(x_1,y_1) = (x+y,y)$.
To be precise, it has the form
\[
x_1=x+y+\Or(y^2), \qquad y_1=y+\Or(y^3).
\]
Since the $n$-th power map is real analytic,
it can be extended to a complex domain of the form
\[
D_{a_*,b_*} :=
\left\{ (x,y) \in \Cset/\Zset \times \Cset : |\Im x|<a_*, |y|<b_* \right\}.
\]
The quantity $a_*$ plays a more important role than $b_*$.
To be precise,
we have proved that, once fixed any $\alpha \in (0,a_*)$ and $L \ge 1$,
there exist $K, q_*>0$ such that
\[
\Delta^{(p,q)}\le K \exp\left(-\frac{2\pi \alpha q}{|np-mq|}\right),
\]
for any relatively prime integers $p$ and $q$ such that
$1 \le |np-mq|\le L$ and $q \ge q_*$.
See Theorem~\ref{thm:Asymptotic} for a more detailed statement.
This upper bound is optimal because $\alpha \in (0,a_*)$.
That is, the exponent $\alpha$ can be taken as close to the
analyticity strip width $a_*$ as desired.
The constant $K$ may explode when $\alpha$ tends to $a_*$,
so, in general, we can not take $\alpha = a_*$.
A similar optimal exponentially small upper bound was obtained
in~\cite{FontichSimo1990} in the setting of the splitting of separatrices
of weakly hyperbolic fixed points of analytical area-preserving maps.
The proof of this optimal bound adds some extra technicalities,
but we feel that the effort is worth it.

The proof is based on two facts.
First, we write the $n$-th power of the exact twist map as the integrable twist map
$(x_1,y_1) = (x+y,y)$ plus an exponentially small remainder on the distance
to the RIC.
See Theorem~\ref{thm:Remainder}.
The size of the remainder is reduced by performing a finite sequence
of changes of variables, but the number of such changes increases
when we approach to the resonant RIC.
This is a classical Neishtadt-like argument~\cite{Neishtadt1981}.
Second, we apply the MacKay-Meiss-Percival action
principle~\cite{MackayMeissPercival1984}, in which the difference
of actions of $(p,q)$-periodic actions is interpreted as an area
on the phase space.

The structure of the paper is the following.
Section~\ref{sec:MainTheorems} is devoted to
state our results in the general context of
analytic exact twist maps.
In Sections~\ref{sec:LengthSpectrum} and~\ref{sec:AreaSpectrum}, we present
the different billiard maps and show how the results in
Section~\ref{sec:MainTheorems} apply.
Sections~\ref{sec:ProofRemainder} and~\ref{sec:ProofThmAsymptotic}
contain the technical proofs.

\section{Main theorems}
\label{sec:MainTheorems}

\subsection{Dynamics close to an analytic resonant RIC}

Let us introduce some notions about exact twist maps defined
on open cylinders.
For a general background on these maps,
we refer to~\cite[\S 9]{KatokHasselblatt1995} and~\cite{Meiss1992}.

Let $\Tset=\Rset/\Zset$ and $I=(r_-, r_+)\subset \Rset$,
for some $-\infty \le r_- < r_+ \le +\infty$.
We will use the coordinates $(s,r)$ for both $\Tset \times I$ and
its universal cover $\Rset \times I$. We refer to $s$ as the angular coordinate.
Let $\omega = - \rmd \lambda$ be an \emph{exact symplectic form}
on the open cylinder $\Tset \times I$ such that $\lambda = \nu(r) \rmd s$
and $\omega = \nu'(r) \rmd s \wedge \rmd r$
for some smooth function $\nu:(r_-,r_+) \to \Rset$.
In particular, $\nu'(r) > 0$.
If $h$ is a real-valued smooth function, $\partial_i h$ denotes the derivative
with respect to the $i$-th variable.

\begin{defi}
A smooth diffeomorphism $g: \Tset \times I \to \Tset \times I$ is an
\emph{exact twist map} when it preserves
the exact symplectic form $\omega = - \rmd \lambda$, has zero flux,
and satisfies the classical  \emph{twist condition} $\partial_2 s_1 (s,r) > 0$,
where $G(s,r)=(s_1,r_1)$ is a lift of $g$.
\end{defi}

Henceforth, the exact symplectic form
$\omega = -\rmd \lambda = \nu'(r) \rmd s \wedge \rmd r$ and
the lift $G$ remain fixed.
We will assume that $\nu(r)$ is analytic when we deal with analytic maps.
If $g$ preserves $\omega$ and $t = \nu(r)$,
then $\lambda = t \rmd s$, $\omega = \rmd s \wedge \rmd t$, and
$g$ preserves the canonical area in the \emph{global Darboux coordinates} $(s,t)$.
It is worth to remark that certain billiard maps can be analytically extended to
the boundaries of their phase spaces in $(s,r)$ variables,
but not in global Darboux coordinates.
For this reason we consider the coordinates $(s,r)$ and
the above exact symplectic forms.
In fact, we could deal with any exact symplectic form,
but we do not need it for the problems we have in mind.

\begin{defi}\label{def:IntersectionProperty}
Let $g:\Tset \times I \to \Tset \times I$ be a continuous map.
The map $g$ has the \emph{intersection property} on $\Tset \times I$
if the image of any closed homotopically non trivial loop of the cylinder
$\Tset \times I$ intersects the loop.
\end{defi}

The intersection property is preserved under
global changes of variables.

\begin{defi}
A \emph{rotational invariant curve (RIC)} of $g$ is a closed loop
$C \subset \Tset \times I$ homotopically non trivial such that $g(C) = C$.
Let $C$ be a RIC of $g$.
Let $m$ and $n$ be two relatively prime integers such that $n\ge 1$.
We say that $C$ is $(m,n)$-\emph{resonant} when $G^n(s,r) = (s+m,r)$
for all $(s,r) \in C$,
and we say that $C$ is \emph{analytic} when
$C = \Graph \gamma := \{ (s,\gamma(s)) : s \in \Tset \}$ for some
analytic function $\gamma: \Tset \to I$.
\end{defi}

If $g: \Tset \times I \to \Tset \times I$ is a diffeomorphism
preserving the exact symplectic form $\omega = -\rmd \lambda$ and has a RIC,
both the zero flux condition and the intersection property
are automatically satisfied.

Let us study the dynamics of an analytic exact twist map
in a neighbourhood of an analytic $(m,n)$-resonant RIC.
First, we note that all points on a $(m,n)$-resonant RIC of $g$
remain fixed under the power map $f = g^n$.
Second, we adapt a classical lemma that appears in several papers
about billiards~\cite{Lazutkin1973,Tabachnikov1995} to our setting.

\begin{lem}\label{lem:Order2}
If $g:\Tset \times I \to \Tset \times I$ is an analytic exact twist map
with an analytic $(m,n)$-resonant RIC $C \subset \Tset \times I$,
then there exist an analytic strip width $a_* > 0$,
an analytic radius $b_* > 0$,
and some analytic coordinates $(x,y)$ such that $C \equiv \{ y=0 \}$
and the power map $f = g^n$ satisfies the following properties:
\begin{enumerate}
\item
It is real analytic on $\Tset \times (-b_*,b_*)$ and
can be analytically extended to the complex domain
\begin{equation}\label{eq:Dab}
D_{a_*,b_*} =
\left\{ (x,y) \in (\Cset/\Zset) \times \Cset : |\Im x| < a_*,\ |y| < b_*
\right\};
\end{equation}
\item
It has the intersection property on the cylinder
$\Tset \times (-b_*,b_*)$; and
\item
It has the form $(x_1,y_1) = f(x,y)$, with
\begin{equation}\label{eq:Order2}
x_1 = x + y + \Or(y^2), \qquad
y_1 = y + \Or(y^3).
\end{equation}
\end{enumerate}
\end{lem}

\proof
If $C = \Graph \gamma$ and $v = r - \gamma(s)$,
then $C \equiv \{v=0\}$ and $(s_1,v_1) = f(s,v)$, with
\begin{equation}\label{eq:varphi_psi}
s_1 = s + \varphi(s)v + \Or(v^2),\qquad
v_1 = v + \psi(s)v^2 + \Or(v^3),
\end{equation}
for some real analytic $1$-periodic functions $\varphi(s)$ and $\psi(s)$.
The twist condition on the RIC implies that $\varphi(s)$ is positive,
since any power of a twist map is locally twist on its smooth
RICs~\cite[Lemma 2.1]{PintodeCarvalhoRamirezRos2013}.
The preservation of $\omega$ implies that
$2 \mu \psi = - (\mu \varphi )'$, where $\mu(s) = \nu'(\gamma(s)) > 0$.
Next, we consider the analytic coordinates $(x,y)$ defined by
\[
x = k \int_0^s \sqrt{\frac{\mu(s)}{\varphi(s)}}\rmd t,\qquad
y = k \sqrt{\mu(s) \varphi(s)} v,\qquad
k^{-1} = \int_0^1 \sqrt{\frac{\mu(t)}{\varphi(t)}} \rmd t.
\]
The constant $k$ has been determined in such a way that
the new angular coordinate $x$ is defined modulus one: $x \in \Tset$.
Clearly, $C \equiv \{ y=0 \}$.
Thus, the coordinates $(x,y)$ cover an open set containing
$\Tset \times \{ 0 \}$, since they are defined in a neighbourhood of $C$.
In particular, $f$ can be analytically extended to
the complex domain $D_{a_*,b_*}$ for some $a_*,b_* > 0$.
Besides,
$f$ has the intersection property on $\Tset \times (-b_*,b_*)$
because the integral of the area form $\omega$ over the region
enclosed between the RIC $C$ and any closed homotopically non
trivial loop should be preserved.
Finally, a straightforward computation shows that $f$ has the form~(\ref{eq:Order2})
in the coordinates $(x,y)$.
\qed

Lemma~\ref{lem:Order2} has, at a first glance, a narrow scope of application
because resonant RICs are destroyed under generic perturbations.
However, the boundaries of the cylinder can be considered
resonant RICs of the extended twist map in many interesting examples.

Let us precise this idea.

\begin{defi}
Let $g: \Tset \times I \to \Tset \times I$ be a continuous map.
If $r_-$ is finite, we say that $C_- = \Tset \times \{r_- \}$ is a
\emph{rigid rotation boundary} when $g$ can be continuously extended to
$\Tset \times [r_-,r_+)$ and its extended lift satisfies
that $G(s,r_-) = s + \varrho_-$ for some \emph{boundary rotation number}
$\varrho_- \in \Rset$,
and we say that $C_-$ is a $(m,n)$-\emph{resonant boundary}
when $\varrho_- = m/n$.
\end{defi}

\begin{defi}
Let $g: \Tset \times I \to \Tset \times I$ be a smooth diffeomorphism.
We say that the \emph{twist condition holds on the boundary} $C_-$
when $g$ can be smoothly extended to $\Tset \times [r_-,r_+)$
and its extended lift satisfies that $\partial_2 s_1(s,r_-) > 0$
for all $s \in \Tset$.
\end{defi}

Analogous definitions can be written for the upper boundary
$C_+ = \Tset \times \{ r_+ \}$ when $r_+$ is finite.
We recall that $G(s,r)=(s_1,r_1)$ is a lift of $g$.
Next, we present a version of Lemma~\ref{lem:Order2} for
the boundaries of the cylinder.
The only remarkable difference is that the exact symplectic form
$\omega = -\rmd \lambda = \nu'(r) \rmd s \wedge \rmd r$ may vanish on the boundaries.

\begin{lem}\label{lem:Order2_Boundaries}
If $g:\Tset \times I \to \Tset \times I$ is an analytic exact twist map
such that $r_-$ is finite,
$C_- = \Tset \times \{ r_- \}$ is a $(m,n)$-resonant boundary,
$g$ can be analytically extended to $\Tset \times [r_-,r_+)$,
$\nu(r)$ can be analytically extended to $[r_-,r_+)$,
and the twist condition holds on $C_-$,
then there exist an analytic strip width $a_* > 0$,
an analytic radius $b_* > 0$,
and some analytic coordinates $(x,y)$ such that $C_- \equiv \{ y=0 \}$
and the power map $f = g^n$ satisfies the properties (i)--(iii) given
in Lemma~\ref{lem:Order2}.

An analogous result holds for the upper boundary $C_+ = \Tset \times \{ r_+ \}$
when $r_+$ is finite.
\end{lem}

\proof
If $v = r - r_-$, then $C_- \equiv \{ v = 0 \}$ and the extended power map
$(s_1,v_1) = f(s,v)$ has the form~(\ref{eq:varphi_psi}) for
some real analytic $1$-periodic functions $\varphi(s)$ and $\psi(s)$.
The twist condition on $C_-$ implies that $\varphi(s)$ is positive.
Since $\nu(r)$ is analytic in $[r_-,r_+)$ and $\nu'(r)$ is positive in $(r_-,r_+)$,
we deduce that $\nu'(r) = \nu_* v^j + \Or(v^{j+1})$
for some $\nu_* > 0$ and some integer $j \ge 0$.
The preservation of $\omega$ implies that $\psi = -\varphi'/(j+2)$.
Next, we consider the analytic coordinates $(x,y)$ defined by
\begin{equation}\label{eq:AnalyticCoordinates_xy}
x = k \int_0^s \varphi^{-(j+1)/(j+2)}(t) \rmd t,\qquad
y = k \varphi^{1/(j+2)}(s) v,
\end{equation}
where the constant $k$ is determined in such a way that the angular
coordinate $x$ is defined modulus one.
The rest of the proof follows the same lines as in Lemma~\ref{lem:Order2}.
We just note that the intersection property on a (real) neighbourhood
of the boundary $C_-$ follows by analytic extension, since $f$ preserves the
exact symplectic form $\omega$ for negative values of $v$ too.
\qed

We unify the resonant RICs studied in Lemma~\ref{lem:Order2}
and the resonant boundaries studied in Lemma~\ref{lem:Order2_Boundaries}
as a single object for the sake of brevity.

\begin{defi}\label{defi:C}
Let $g: \Tset \times I \to \Tset \times I$ be an analytic map.
Let $m$ and $n \ge 1$ be two relatively prime integers.
Let $a_* >0$ and $b_* > 0$.
We say that $C \subset \Tset \times I$
(respectively, $C = \Tset \times \{ r_- \}$ for a finite $r_-$
            or $C = \Tset \times \{ r_+ \}$ for a finite $r_+$)
is a \emph{$(a_*,b_*)$-analytic $(m,n)$-resonant RIC}
when $C$ is an analytic $(m,n)$-resonant RIC (respectively, $(m,n)$-resonant boundary)
and there exist some coordinates $(x,y)$ such that $C \equiv \{ y = 0 \}$ and
the power map $f = g^n$ satisfies the properties (i)--(iii) stated
in Lemma~\ref{lem:Order2}.
\end{defi}

The map~(\ref{eq:Order2}) can be viewed as a perturbation of
the integrable twist map
\begin{equation}\label{eq:IntegrableMap}
x_1 = x + y, \qquad y_1 = y.
\end{equation}
We want to reduce the size of the nonintegrable terms
$\Or(y^2)$ and $\Or(y^3)$ as much as possible.
We can reduce them through normal form steps up to
any desired order; see Lemma~\ref{lem:BeyondAllOrders}.
Thus, the nonintegrable part of the dynamics is beyond all order in $y$;
that is, in the distance to the resonant RIC (or resonant boundary).
A general principle in conservative dynamical systems states that
beyond all order phenomena are often exponentially small in the
analytic category.
Our goal is to write the map as an exponentially small
perturbation in $y$ of the integrable twist map~(\ref{eq:IntegrableMap}).
The final result is stated in the following theorem.

\begin{thm}\label{thm:Remainder}
Let $f:\Tset \times I \to \Tset \times I$ be an analytic map with
a $(a_*,b_*)$-analytic RIC $C$ of fixed points such that
$C \subset \Tset \times I$,
$C = \Tset \times \{ r_- \}$ for a finite $r_-$,
or $C = \Tset \times \{ r_+ \}$ for a finite $r_+$.
Let $m \ge 2$ be an arbitrary order.
Let $\alpha \in (0,a_*)$.
There exist constants $K > 0$ and $b'_* \in (0,b_*)$ such that,
if $b \in (0, b'_*)$, then there exists an analytic change of variables
$(x,y) = \Phi(\xi,\eta)$ such that:
\begin{enumerate}
\item
It is uniformly (with respect to $b$) close to the identity on
$\Tset \times (-b,b)$. That is,
\[
x = \xi + \Or(\eta),\qquad
y = \eta + \Or(\eta^2),\qquad
\det[\rmD \Phi(\xi,\eta)] = 1 +\Or(\eta),
\]
for all $(\xi,\eta) \in \Tset \times (-b,b)$,
where the $\Or(\eta)$ and $\Or(\eta^2)$ terms are uniform in $b$; and
\item
The transformed map $(\xi,\eta) \mapsto (\xi_1,\eta_1)$ is real analytic
and has the intersection property on the cylinder $\Tset \times (-b,b)$.
Besides, it has the form
\begin{equation}\label{eq:NeishtadtForm}
\xi_1  = \xi + \eta + \eta^m g_1(\xi,\eta), \qquad
\eta_1 = \eta + \eta^{m+1} g_2(\xi,\eta),
\end{equation}
where $|g_j(\xi,\eta)|\le K\rme^{-2\pi \alpha / b}$
and $|\partial_i g_j(\xi,\eta)| \le K b^{-2}$
for all $(\xi,\eta)\in \Tset \times (-b,b)$.
\end{enumerate}
\end{thm}

The proof can be found in Section~\ref{sec:ProofRemainder}.

Consider a perturbed Hamiltonian system which is close to an integrable system.
It is known that, under the appropriate nondegeneracy conditions,
the measure of the set of tori which decompose under the perturbation
can be bounded from above by a quantity of order $\sqrt{\epsilon}$,
$\epsilon$ being the perturbation parameter~\cite{Neishtadt1981,Poschel1982}.
Neishtadt~\cite{Neishtadt1981} also considered a context where the perturbation
becomes exponentially small in some parameter $\epsilon$ and hence
the measure of the complementary set
which is cut out from phase space by the invariant tori is of order
$\rme^{-c/\epsilon}$, $c$ being a positive constant.
This argument could be applied to our context.
First, any neighbourhood of an analytic resonant RIC (or resonant boundary)
of an analytic exact twist map contains infinitely many RICs.
Second, the area of the complementary of the RICs in any of such neighbourhoods
is exponentially small in the size of the neighbourhood.
Third, the gaps between the RICs are exponentially small in their distance
to the resonant RIC.
The first result follows from the classical Moser twist theorem~\cite{SiegelMoser1995}.
The others follow from the ideas explained above.

\subsection{Difference of periodic actions}
\label{ssec:DifferenceOfPeriodicActions}

Let $\omega = -\rmd \lambda$, with $\lambda = \nu(r) \rmd s$,
be a fixed exact symplectic form on the open cylinder
$\Tset \times I$, with $I=(r_-,r_+)$.
Let $g:\Tset\times I \to \Tset\times I$ be an exact twist map
that can be extended as rigid rotations of angles $\varrho_-$
and $\varrho_+$ to the boundaries $C_- = \Tset \times \{ r_- \}$ and
 $C_+ = \Tset \times \{ r_+ \}$, respectively.
We know that $\varrho_- < \varrho_+$ from the twist condition.
Let $G:\Rset \times I \to \Rset \times I$ be a fixed lift of $g$.
Let $E=\{(s,s_1)\in \Rset^2 : \varrho_- < s_1 - s < \varrho_+ \}$.
Then there exists a function $h: E \to \Rset$,
determined modulo an additive constant, such that
\[
G^* \lambda - \lambda = \rmd h.
\]
The function $h$ is called the \emph{Lagrangian} or
\emph{generating function} of $g$.

Let $p$ and $q$ be two relatively prime integers
such that $\varrho_-<p/q<\varrho_+$ and $q\ge 1$.
A point $(s,r) \in \Rset \times I$ is $(p,q)$\emph{-periodic}
when $G^q(s,r)=(s+p,r)$.
The corresponding point $(s,r) \in \Tset \times I$
is a periodic point of period $q$ by $g$ that is translated
$p$ units in the base by the lift.
A $(p,q)$-periodic orbit is \emph{Birkhoff} when it is ordered
around the cylinder in the same way that the orbits of the
rigid rotation of angle $p/q$.
See~\cite{KatokHasselblatt1995} for details.
The Poincar\'e-Birkhoff Theorem states that
there exist at least two different Birkhoff $(p,q)$-periodic
orbits~\cite{KatokHasselblatt1995,Meiss1992}.

Let $O=\{(s_k, r_k)\}_{k\in\Zset}$ be a $(p,q)$-periodic orbit.
Its $(p,q)$\emph{-periodic action} is
\[
W^{(p,q)}[O] = h(s_0,s_1) + h(s_1,s_2) + \cdots + h(s_{q-1},s_0 + p).
\]
Our goal is to establish an exponentially small bound for
the non-negative quantity
\[
\Delta^{(p,q)} =
\sup_{O,\bar{O}\in\mathcal{O}^{(p,q)}_g}
\left| W^{(p,q)}[\bar{O}] - W^{(p,q)}[O] \right|,
\]
where $\mathcal{O}^{(p,q)}_g$ denotes the set of all Birkhoff
$(p,q)$-periodic orbits of the exact twist map $g$.
The difference of $(p,q)$-periodic actions can be interpreted as
the $\omega$-area of certain domains.

Let us explain it.

Let $O=\{(s_k,r_k)\}_{k\in\Zset}$ and $\bar{O}=\{(\bar{s}_k,\bar{r}_k)\}_{k\in\Zset}$ be
two $(p,q)$-periodic orbits.
We can assume, without loss of generality, that  $0<\bar{s}_0-s_0<1$.
Let $L_0$ be a curve from $(s_0,r_0)$ to $(\bar{s}_0,\bar{r}_0)$
contained in $\Rset \times I$.
Set $L_k=G^k(L_0)$.
The curves $L_0$ and $L_q$ have the same endpoints in $\Tset\times I$.
Let us assume that these two curves have no topological crossing on the
cylinder $\Tset \times I$ and let $B \subset \Tset \times I$ be the
domain enclosed between them.

Observe that
$ \int_{L_{k+1}} \lambda- \int_{L_k} \lambda =
\int_{L_k} (G^* \lambda - \lambda) =
\int_{L_k} \rmd h =
 h(\bar{s}_k, \bar{s}_{k+1}) - h(s_k, s_{k+1})$.
Hence, $ \sum_{k=0}^{q-1} \big( h(\bar{s}_k, \bar{s}_{k+1}) - h(s_k, s_{k+1}) \big) =
\int_{L_q} \lambda - \int_{L_0} \lambda = \int_{L_q - L_0} \lambda = \pm \int_B \omega$,
where the sign $\pm$ depends on the orientation of the closed path $g^q(L_0) - L_0$,
but we do not need it, because we take absolute values in both sides of the
previous relation:
\begin{equation}\label{eq:MMP}
\left| W^{(p,q)}[\bar{O}] - W^{(p,q)}[O] \right| =
\int_B \omega =: \Area_\omega [B].
\end{equation}
These arguments go back to the
\emph{MacKay-Meiss-Percival action principle}~\cite{MackayMeissPercival1984, Meiss1992}.
If the curves $L_0$ and $g^q(L_0)$ have some topological crossing,
then the domain $B$ has several connected components,
in which case $\left| W^{(p,q)}[\bar{O}] - W^{(p,q)}[O] \right| \le
\int_B \omega =: \Area_\omega [B]$, because the sign in front of
the integral of the area form $\omega$ depends on the connected component.

If the analytic exact twist map $g$ has a $(m,n)$-resonant RIC,
then
\[
\Delta^{(m,n)} = 0.
\]
Indeed, we can take a segment of the RIC as the curve $L_0$ used in
the previous construction in such a way that $g^n(L_0) = L_0$ and
$\Area_\omega[B] = 0$.
It turns out that the differences of $(p,q)$-periodic actions of $g$
are exponentially small when $p/q$ is ``sufficiently close'' to $m/n$.
The meaning of ``sufficiently close'' is clarified in the following theorem.
See also Remark~\ref{rem:SufficientlyClose}.

\begin{thm}\label{thm:Asymptotic}
Let $g:\Tset\times I\to\Tset\times I$ be an analytic exact twist map
that can be extended as rigid rotations of angles $\varrho_-$
and $\varrho_+$ to the boundaries $C_- = \Tset \times \{ r_- \}$ and
 $C_+ = \Tset \times \{ r_+ \}$, respectively.
We also assume that $g$ has a
$(a_*,b_*)$-analytic $(m,n)$-resonant RIC $C \subset \Tset \times I$,
$C = \Tset \times \{ r_- \}$ for a finite $r_-$,
or $C = \Tset \times \{ r_+ \}$ for a finite $r_+$.
Let $\alpha \in (0,a_*)$ and $L \ge 1$.
There exist $K, q_* > 0$ such that
\begin{equation}\label{eq:DeltaBound}
\Delta^{(p,q)}\le K \exp\left(-\frac{2\pi \alpha q}{|np-mq|}\right),
\end{equation}
for all relatively prime integers $p$ and $q$ with $1 \le |np-mq|\le L$
and $q\ge q_*$.
\end{thm}

The proof has been placed at Section~\ref{sec:ProofThmAsymptotic}.

\begin{remark}\label{rem:TendingToBoundaryRics}
If $C = \Tset \times \{ r_- \}$ (respectively, $C = \Tset \times \{ r_+ \}$),
the bound~(\ref{eq:DeltaBound}) only makes sense as $p/q$ tends to the
boundary rotation number $\varrho_- = m/n$ (respectively, $\varrho_+ = m/n$)
from the right: $p/q \to (m/n)^+$
(respectively, from the left: $p/q \to (m/n)^-$).
We will see several examples of this situation in
Sections~\ref{sec:LengthSpectrum} and~\ref{sec:AreaSpectrum}.
\end{remark}

\begin{remark}\label{rem:InfiniteValues}
Infinite values of $r_\pm$ can also be dealt with.
For instance, if $r_+ = +\infty$, then we consider the coordinate $v = 1/r$
and we assume that:
1) $C_+ = \{ v=0 \}$ is a $(m,n)$-resonant boundary;
2) $g$ can be analytically extended to $v \ge 0$;
3) $\omega = \beta(v) \rmd s \wedge \rmd v$ with
   $\beta(v) = \beta_* v^j + \Or(v^{j+1})$ for some real number $\beta_* \neq 0$
   and some integer $j \neq -2$;
4) $f = g^n$ has the form~(\ref{eq:varphi_psi}) for some real analytic $1$-periodic
   functions $\varphi(s)$ and $\psi(s)$; and
5) $\varphi(s)$ is negative.
Then the change $(s,r) \mapsto (s,v=1/r)$ transforms the infinite setup
into a finite one, so Theorem~\ref{thm:Asymptotic} still holds.
We just mention two subtle points.
First, the change~(\ref{eq:AnalyticCoordinates_xy})
is not well-defined when $j = -2$.
Thus, $\beta(v)$ can have a pole at $v = 0$, provided it is not a double one.
Second, $\varphi(s)$ is negative because the change $(s,r) \mapsto (s,v=1/r)$
reverses the orientation, and so, it changes the sign of the twist.
\end{remark}

\begin{remark}\label{rem:SufficientlyClose}
Condition $|np-mq| \le L$ implies that $|p/q-m/n| = \Or(1/q)$
as $q \to +\infty$.
\end{remark}

\begin{remark}
In many applications, $C$ is a RIC of fixed points; that is, a $(0,1)$-resonant RIC.
Then Theorem~\ref{thm:Asymptotic} implies that $\Delta^{(p,q)}$ is exponentially small
in the period $q$ when $p$ remains uniformly bounded.
To be precise, if $\alpha \in (0,a_*)$ and $L \ge 1$,
there exist $K, q_* > 0$ such that
\[
\Delta^{(p,q)} \le K \rme^{-2\pi \alpha q/|p|}
\]
for all relatively prime integers $p$ and $q$ with
$q \ge q_*$ and $0 \neq |p| \le L$.
\end{remark}

\section{On the length spectrum of analytic convex domains}
\label{sec:LengthSpectrum}

\subsection{Convex billiards}

We recall some well-known results about billiards that can be found
in~\cite{KozlovTreschev1991,Tabachnikov1995_Billiards,KatokHasselblatt1995}.

Let $\Gamma$ be a closed strictly convex curve in the Euclidean plane $\Rset^2$.
We assume, without loss of generality, that $\Gamma$ has length one.
Let $s \in \Tset$ be an arc-length parameter on $\Gamma$.
Set $I = (r_-,r_+) = (0,\pi)$.
Let $f:\Tset \times I \to \Tset \times I$, $(s_1,r_1) = f(s,r)$,
be the map that models the billiard dynamics inside $\Gamma$ using
the \emph{Birkhoff coordinates} $(s,r)$,
where $s \in \Tset$ determines the impact point on the curve,
and $r \in I$ denotes the angle of incidence.

The map $f$ preserves the exact symplectic form
$\omega = \nu'(r) \rmd s \wedge \rmd r = \sin r  \rmd s \wedge \rmd r$
and has the intersection property.
Indeed, $f: \Tset \times I \to \Tset \times I$ is an exact twist map
with boundary rotation numbers $\varrho_- = 0$ and $\varrho_+ = 1$.
Besides, its Lagrangian is the distance between consecutive impact points.
Finally, $f$ is analytic when $\Gamma$ is analytic.

Any $(p,q)$-periodic orbit on the billiard map forms a
a closed inscribed polygon with $q$ sides that makes $p$ turns inside $\Gamma$.
Since the Lagrangian of the billiard map is the distance between consecutive
impact points,
the periodic action of a periodic orbit is just the total
length of the corresponding polygon.
Therefore, the supremum action difference among $(p,q)$-periodic billiard orbits
is the supremum length difference among inscribed billiard $(p,q)$-polygons.

\subsection{Study close to the boundary of the billiard table}

Let us check that the boundary $C_- = \Tset \times \{ 0 \}$
satisfies the hypotheses stated in Theorem~\ref{thm:Asymptotic}.

\begin{pro}\label{pro:ClassicalExtensionToBoundary}
If $f$ is the billiard map associated to an analytic strictly convex curve
$\Gamma$, the boundary $C_-$ is a $(a_*,b_*)$-analytic $(0,1)$-resonant
RIC of $f$ for some $a_* > 0$ and $b_* > 0$.
\end{pro}

\proof
The lower boundary $C_-$ is a $(0,1)$-resonant because
$\varrho_- = 0$.
Clearly, $r_- = 0$ is finite and the function $\nu(r) = \sin r$
can be analytically extended to $[0,\pi)$.
Hence, in order to end the proof we only need to prove that the
billiard map $f$ can be analytically extended to $\Tset \times [0,\pi)$
and the twist condition holds on $C_-$; see Lemma~\ref{lem:Order2_Boundaries}.

Let $m: \Tset \to \Gamma$, $\kappa: \Tset \to (0,+\infty)$, and
$n: \Tset \to \Rset^2$ be an arc-length parametrization,
the curvature, and a unit normal vector of the
analytic strictly convex curve linked by the relation
$m''(s) = \kappa(s) n(s)$.

We write the angle of incidence as a function of
consecutive impact points: $r = r(s,s_1)$.
If $s_1 \neq s$, then $r(s,s_1)$ is analytic,
$r(s,s_1) \in (0,\pi)$, and $\partial_2 r(s,s_1) > 0$.
The last property follows from the twist character of the billiard map.
Let us study what happens on $C_-$;
or, equivalently, what happens in the limit $s_1 \to s^+$.
We note that $r$ is the angle between $m'(s)$ and $m(s_1)-m(s)$.
Hence,
\begin{equation}
\fl
\tan r =
\frac{\sin r}{\cos r} =
\frac{\det\left(m'(s),m(s_1)-m(s)\right)}
     {\langle m'(s),m(s_1)-m(s) \rangle} =
\frac{\int_0^1 \det\left(m'(s),m'(s+t(s_1-s))\right)\rmd t}
     {\int_0^1\langle m'(s),m'(s+t(s_1-s))\rangle\rmd t}.
\end{equation}
This expression shows that $r(s,s_1)$ is analytic when $s_1 = s$ and $r(s,s) = 0$.
Besides,
\begin{equation}\label{eq:TwistBoundaryBilliard}
\lim_{s_1 \to s^+} \partial_2 r(s,s_1) =
\lim_{s_1 \to s^+} \cos^2 r(s,s_1) \partial_2 \tan r(s,s_1) =
\kappa(s)/2 > 0,
\end{equation}
for all $s \in \Tset$, which implies that $(s,s_1) \mapsto (s,r)$
is an analytic diffeomorphism that maps a neighbourhood of the diagonal in $\Tset^2$
to a neighborhood of the lower boundary $C_- = \Tset \times \{ 0 \}$.
Next, we write the billiard map $(s_1,r_1) = f(s,r)$ as
the composition of three maps:
\[
(s,r) \mapsto (s,s_1) \mapsto (s_1,s) \mapsto (s_1,r_1).
\]
We have already seen that the first map is analytic in a neighborhood of $C_-$.
Clearly, the second map is analytic.
The third one is also analytic,
because $r_1 = r_1(s,s_1) = \pi-r(s_1,s)$
by definition of billiard map.
Thus, $f$ can be analytically extended to $\Tset \times [0,\pi)$.
The twist condition on $C_-$ follows from
inequality~(\ref{eq:TwistBoundaryBilliard}).
\qed

We can compare the coordinates~(\ref{eq:AnalyticCoordinates_xy})
with the coordinates defined by Lazutkin in~\cite{Lazutkin1973}.
We note that $\nu'(r) = \sin r = r + \Or(r^2)$, so $j = 1$
in the change~(\ref{eq:AnalyticCoordinates_xy}).
If $\rho(s) = 1/\kappa(s)$ is the radius of curvature of $\Gamma$,
then the Taylor expansion around $r = 0$ of the billiard map is
\[
\left\{
\begin{array}{l}
s_1  =  s + 2\rho(s)r + 4\rho(s)\rho'(s)r^2/3 + \Or(r^3), \\
r_1  =  r - 2\rho'(s)r^2/3 +
\left(4(\rho'(s))^2/9 - 2\rho(s)\rho''(s)/3\right)r^3 + \Or(r^4).
\end{array}
\right.
\]
From this Taylor expansion,
Lazutkin deduced that the billiard map takes the form
\[
x_1 = x + y + \Or(y^3),\qquad y_1 = y + \Or(y^4)
\]
in the analytic \emph{Lazutkin coordinates} $(x,y)$ defined by
\[
x = k \int_0^s \rho^{-2/3}(t) \rmd t,\qquad
y = 4 k \rho^{1/3}(s)\sin(r/2),\quad
k^{-1} = \int_0^1 \rho^{-2/3}(t)\rmd t.
\]
The constant $k$ has, once more, been determined in such a way that
the new angular coordinate $x$ is defined modulus one.
Lazutkin's results are more refined, he wrote the billiard map
as a smaller perturbation of the integrable twist map
$(x_1,y_1) = (x+y,y)$, but we do not need it.

By direct application of Proposition~\ref{pro:ClassicalExtensionToBoundary}
and Theorem~\ref{thm:Asymptotic}, we get the following exponentially small upper bound of
the quantities $\mathcal{L}^{(p,q)}$ defined in the introduction.

\begin{thm}\label{thm:ClassicalBoundary}
Let $\Gamma$ be an analytic strictly convex curve in the Euclidean plane.
Let $a_* > 0$ be the analyticity strip width of the lower boundary $C_-$.
Let $\alpha \in (0,a_*)$ and $L \ge 1$.
There exist a constant $K > 0$ and a period $q_* \ge 1$ such that
\[
\mathcal{L}^{(p,q)} \le K \rme^{-2\pi \alpha q / p},
\]
for all relatively prime integers $p$ and $q$ with $q \ge q_*$ and $0 < p \le L$.
\end{thm}

The same exponentially small upper bound holds for analytic
geodesically strictly convex curves on surfaces of constant curvature,
where the billiard trajectories are just broken geodesics.
Billiard maps on the Klein model of the hyperbolic plane $\Hset^2$
and on the positive hemisphere $\Sset^2_+$ have been studied,
for instance, in~\cite{Coutinho},
where it is shown that they are exact twist maps
with the same boundary rotation numbers as in the Euclidean case.
Therefore, by local isometry arguments,
we can write a version of Theorem~\ref{thm:ClassicalBoundary} on any
surface of constant curvature.

\subsection{Billiard tables of constant width}\label{sec:constantwidth}

\begin{defi}
A smooth closed convex curve is of \emph{constant width}
if and only if it has a chord in any direction perpendicular
to the curve at both ends.
\end{defi}

Billiards inside convex curves of constant width have
a nice property~\cite{Knill1998,Gutkin2012}.
Let us explain it.

The billiard map associated to a smooth convex curve of constant width has
the horizontal line $\Tset \times \{ \pi/2\}$ as a resonant $(1,2)$-RIC.
Any trajectory belonging to that RIC is orthogonal to the curve
at its two endpoints.
Due to the variational formulation, all the $(1,2)$-periodic orbits are extrema
of the $(1,2)$-periodic action and, thus, all $(1,2)$-periodic trajectories
have the same length, which is the reason we refer to them
as constant width curves.

Another characterization of constant width curves is the following.
We reparametrize the curve by using the angle $\varphi \in (0,2\pi)$
between the tangent vector at a point in the curve and some fixed line.
Let $\rho(\varphi)$ be the radius of curvature at this point.
The curve has constant width if and only if the Fourier series
of $\rho(\varphi)$ contains no other even coefficients than
the constant term.
Thus, the space of analytic constant width curves
has infinite dimension and codimension.

Theorem~\ref{thm:Asymptotic} can be directly applied
in this context and we have the following result.

\begin{thm}
Let $g: \Tset \times I \to \Tset \times I$ be the billiard map
of an analytic strictly convex curve of constant width.
Let $a_* > 0$ be the analyticity strip width of
the $(1,2)$-resonant RIC of $g$.
Let $\alpha \in (0,a_*)$ and $L \ge 1$.
There exist a constant $K>0$ and a period $q_* \ge 1$ such that
\[
\mathcal{L}^{(p,q)} \le K\exp\left(-\frac{\pi\alpha}{|p/q-1/2|}\right),
\]
for all relatively prime integers $p$ and $q$ such that
$1 \le |2p-q|\le L$ and $q\ge q_*$.
\end{thm}

One could try to generalize constant width billiards,
where $\Tset \times \{\pi/2 \}$ is a $(1,2)$-resonant RIC,
to \emph{constant angle tables}, where
$\Tset \times \{ r_0 \}$ is assumed to be a $(m,n)$-resonant RIC.
However, the only table such that $\Tset \times \{ r_0 \}$
is a $(m,n)$-resonant RIC, with $(m,n)\neq(1,2)$, is the circle.
See~\cite{Gutkin2012,Cyr2012}.
By the way, Theorem~\ref{thm:Asymptotic} applies to this case but,
since the circular billiard is integrable,
$\mathcal{L}^{(m,n)} \equiv 0$, for all $(m,n)$.
In fact, the circular billiard map is globally conjugated
to the integrable twist map~(\ref{eq:IntegrableMap}).

There are more billiard tables with resonant RICs,
but their RICs are not horizontal.
For instance,
the elliptic table has all possible $(m,n)$-resonant RICs,
but the $(1,2)$-resonant one.
Hence, in this case, $\Delta^{(m,n)}=0$.
Baryshnikov and Zharnitsky~\cite{BaryshnikovZharnitsky2006} proved that an ellipse
can be infinitesimally perturbed so that any chosen resonant RIC will persist.
Innami~\cite{Innami1988} found a condition on the billiard table
that guarantees the existence of a $(1,3)$-resonant RIC.
However, Theorem~\ref{thm:Asymptotic} can not be applied in such cases,
because both the Baryshnikov-Zharnitsky and the Innami constructions
are done in the smooth category,
where we can only claim that $\mathcal{L}^{(p,q)}$ is beyond all order
in the difference between rotation numbers.

\section{On the area spectrum of analytic convex domains}
\label{sec:AreaSpectrum}

\subsection{Dual billiards}\label{ssec:DualBilliards}

We recall some well-known facts about dual billiards that can be found
in~\cite{Boyland1996,GutkinKatok1995,Tabachnikov1995_Billiards}.

Let $\Gamma$ be a strictly convex closed curve in
the Euclidean plane $\Rset^2$.
Let $\mho$ be the unbounded component of $\Rset^2 \setminus \Gamma$.
The dual billiard map $f: \mho \to \mho$ is defined as follows:
$f(z)$ is the reflection of $z$ in the tangency point of the oriented
tangent line to $\Gamma$ through $z$.
This map is area-preserving.
Next, we introduce the \emph{envelope coordinates}
$(\alpha,r) \in \Tset_\star \times I$ of a point $z \in \mho$.
In this section, $\Tset_\star = \Rset / 2\pi \Zset$ and $I = (0,+\infty)$.
We recall that $\Tset = \Rset/\Zset$.

Given a point $z \in \mho$, let $\alpha \in \Tset_\star$ be the angle made by the
positive tangent line to $\Gamma$ in the direction of $z$ with a fixed
direction of the plane,
and let $r \in I$ be the distance along this line from $\Gamma$ to $z$.
See Figure~\ref{fig:DualEnvelopeCoordinates}.

\iffigures
\begin{figure}
\centering
\includegraphics[width=0.8\textwidth]{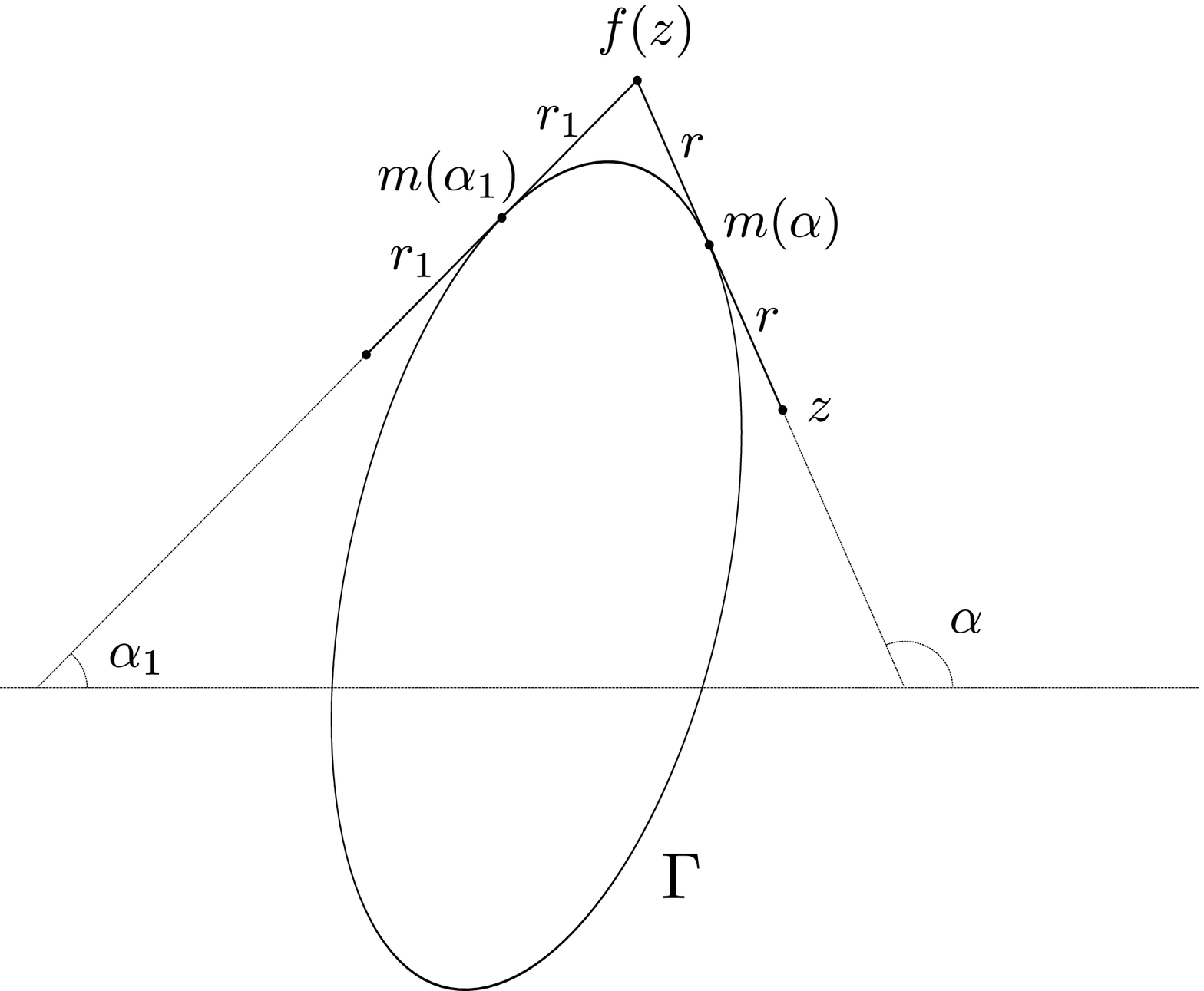}
\caption{The envelope coordinates $(\alpha,r)$ and
         the dual billiard map $f:\mho \to \mho$.}
\label{fig:DualEnvelopeCoordinates}
\end{figure}
\fi

The dual billiard map preserves the exact symplectic form
$\omega = \nu'(r) \rmd \alpha \rmd r = r \rmd \alpha \rmd r$,
and it has the intersection property.
Indeed, $f: \Tset_\star \times I \to \Tset_\star \times I$,
$(\alpha_1,r_1) = f(\alpha,r)$, is an exact twist map with boundary
rotation numbers
$\varrho_- = 0$ and $\varrho_+ = \pi$.
Its Lagrangian is the area enclosed by $\Gamma$
and the tangent lines through the points on $\Gamma$ with
coordinates $\alpha$ and $\alpha_1$.

Any $(p,q)$-periodic orbit on the dual billiard map forms a closed
circumscribed polygon with $q$ sides that makes $p$ turns outside $\Gamma$.
Since the Lagrangian of the dual billiard map is the above-mentioned area,
the periodic action of a periodic orbit is just the area enclosed between
the corresponding polygon and $\Gamma$,
taking into account some multiplicities when $p \ge 2$.
Therefore, the supremum action difference among $(p,q)$-periodic dual billiard orbits
is the supremum area difference among circumscribed dual billiard $(p,q)$-polygons.

\subsection{Study close to the curve}

We note that $r \to 0^+$ when the point $z \in \mho$ approaches
to the curve $\Gamma$.
Therefore, in order to study the dual billiard dynamics close to $\Gamma$,
we must study the dual billiard map $f$ in a neighbourhood of
the lower boundary $C_- = \Tset_\star \times \{ 0 \}$ of $\Tset_\star \times I$.
Let us check that $C_-$ satisfies the hypotheses stated in
Theorem~\ref{thm:Asymptotic}.

\begin{pro}\label{pro:DualExtensionToBoundary}
If $f$ is the dual billiard map associated to an analytic strictly convex curve
$\Gamma$, the lower boundary $C_-$ is an $(a_*,b_*)$-analytic $(0,1)$-resonant
RIC of $f$ for some $a_* > 0$ and $b_* > 0$.
\end{pro}

\proof
The lower boundary $C_-$ is a $(0,1)$-resonant because $\varrho_- = 0$.
Clearly, $r_- = 0$ is finite and the function $\nu(r) = r$
can be analytically extended to $[0,+\infty)$.
Hence, in order to end the proof we only need to prove that the dual
billiard map $f$ can be analytically extended to $\Tset_\star \times [0,+\infty)$
and the twist condition holds on $C_-$; see Lemma~\ref{lem:Order2_Boundaries}.

We write the distance $r$ as a function of consecutive tangent points:
$r = r(\alpha,\alpha_1)$.
We know that $r(\alpha,\alpha_1)$ is analytic, $r(\alpha,\alpha_1) > 0$,
and $\partial_2 r(\alpha,\alpha_1) > 0$
when $\alpha_1 \not \in \{ \alpha,\alpha+\pi \}$.\footnote{Tangent
lines through points with coordinates $\alpha$ and $\alpha + \pi$
are parallel, so $\lim_{\alpha_1 - \alpha \to \pi} r(\alpha,\alpha_1) = +\infty$.}
The last property follows from the twist character of the dual billiard map.

Let us study what happens on $C_-$; or, equivalently,
what happens in the limit $\alpha_1 \to \alpha$.
Let $\rho: \Tset_\star \to (0,+\infty)$ be the radius of curvature of $\Gamma$
in the angular coordinate $\alpha$.
Set $\alpha_1 = \alpha + \delta$.
 From Boyland~\cite{Boyland1996}, we know that
\begin{equation}\label{eq:RfromAlphas}
r(\alpha,\alpha_1) =
\frac{\int_{\alpha}^{\alpha_1}\sin(v - \alpha) \rho(v) \rmd v}
     {\sin (\alpha_1 - \alpha)} =
\frac{\delta }{\sin \delta}
\int_0^1 \sin(\delta t) \rho(\alpha + \delta t) \rmd t.
\end{equation}
This expression shows that $r(\alpha,\alpha_1)$ is analytic when
$\alpha_1 = \alpha$ and $r(\alpha,\alpha) = 0$.
By derivating relation~(\ref{eq:RfromAlphas}) with respect to $\alpha_1$,
we get that
\begin{equation}\label{eq:TwistBoundaryDualBilliard}
\partial_2 r (\alpha,\alpha) = \rho(\alpha)/2 > 0,
\end{equation}
for all $\alpha \in \Tset_\star$,
which implies that $(\alpha,\alpha_1) \mapsto (\alpha,r)$
is an analytic diffeomorphism that maps a neighbourhood of the diagonal
in $\Tset_\star^2$ to a neighborhood of the lower
boundary $C_- = \Tset_\star \times \{ 0 \}$.
Next, we write the billiard map $(\alpha_1,r_1) = f(\alpha,r)$ as
the composition of three analytic maps:
\[
(\alpha,r) \mapsto
(\alpha,\alpha_1) \mapsto
(\alpha_1,\alpha) \mapsto
(\alpha_1,r_1).
\]
Thus, $f$ can be analytically extended to $\Tset_\star \times [0,+\infty)$.
The twist condition on $C_-$ follows from
inequality~(\ref{eq:TwistBoundaryDualBilliard}).
\qed

We can compare the coordinates~(\ref{eq:AnalyticCoordinates_xy})
with the coordinates defined by Tabachnikov in~\cite{Tabachnikov1995}.
We note that $\nu'(r) = r$, so $j = 1$ in the change~(\ref{eq:AnalyticCoordinates_xy}).
If $\kappa(\alpha)$ is the curvature of $\Gamma$,
then the Taylor expansion around $r = 0$ of the dual billiard map is
\[
\alpha_1  =  \alpha + 2\kappa(\alpha) r + \Or(r^2), \qquad
r_1  =  r - 2 \kappa'(\alpha) r^2/3 + \Or(r^3).
\]
(Tabachnikov wrote this Taylor expansion in coordinates $(s,r)$,
 where $s$ is an arc-length parameter, but his result can be easily adapted.)
From this Taylor expansion,
one deduces that the dual billiard map takes the form~(\ref{eq:Order2})
in the analytic \emph{Tabachnikov coordinates}
\[
x = k \int_0^\alpha \kappa^{-2/3}(v) \rmd v,\qquad
y = 2 k \kappa^{1/3}(\alpha) r,\quad
k^{-1} = \int_0^{2\pi} \kappa^{-2/3}(v)\rmd v.
\]
The constant $k$ has been determined in such a way that
$x$ is defined modulus one: $x \in \Tset$.

We get the following exponentially small upper bound of the quantities
$\mathcal{A}^{(p,q)}$ defined in the introduction by direct application of
Proposition~\ref{pro:DualExtensionToBoundary} and Theorem~\ref{thm:Asymptotic}.

\begin{thm}\label{thm:DualBoundary}
Let $\Gamma$ be an analytic strictly convex curve in the Euclidean plane.
Let $a_* > 0$ be the analyticity strip width of the lower boundary $C_-$.
Let $\alpha \in (0,a_*)$ and $L \ge 1$.
There exist a constant $K > 0$ and a period $q_* \ge 1$ such that
\[
\mathcal{A}^{(p,q)} \le K \rme^{-2\pi \alpha q / p},
\]
for all relatively prime integers $p$ and $q$ with $q \ge q_*$ and $0 < p \le L$.
\end{thm}

Tabachnikov~\cite{Tabachnikov2002} studied the dual billiard map in
the hyperbolic plane $\Hset^2$,
and extended the asymptotic expansion~(\ref{eq:AsymptoticArea}) to that new setting.
He also claimed that there exists an analogous formula for dual billiards
on the unit sphere $\Sset^2$.
Therefore, by local isometry arguments,
we can write a version of Theorem~\ref{thm:DualBoundary} on any
surface of constant curvature.

\subsection{Study far away from the curve}

We note that $r \to +\infty$ when the point $z \in \mho$ moves
away from the curve $\Gamma$.
We use the coordinates $(\alpha,v)$ to work at infinity,
where $v = 1/r$ and $(\alpha,r) \in \Tset_\star \times I$
are the coordinates introduced in Subsection~\ref{ssec:DualBilliards}.
The exact symplectic form $\omega = r \rmd \alpha \rmd r$ becomes
$\omega = -v^{-3} \rmd \alpha \rmd v$ in coordinates $(\alpha,v)$.
Tabachnikov~\cite{Tabachnikov1995} realized that the dual billiard map
at infinity can be seen as a map defined in a neighbourhood
of the $(1,2)$-resonant RIC $C_+ \equiv \{ v = 0 \}$.
To be precise, he saw that the dual billiard map can be
analytically extended to $v \ge 0$, its square has the form
\[
\alpha_1 = \alpha + \varphi(\alpha) v + \Or(v^2),\qquad
v_1 = v + \psi(\alpha) v^2 + \Or(v^3),
\]
for some real analytic $1$-periodic functions $\varphi(\alpha)$ and $\psi(\alpha)$,
and $\varphi(\alpha)$ is negative.
Hence, the following result is deduced from Remark~\ref{rem:InfiniteValues}.

\begin{thm}\label{thm:DualInfinity}
Let $\Gamma$ be an analytic strictly convex curve in the Euclidean plane.
Let $a_* > 0$ be the analyticity strip width of the boundary $C_+$.
Let $\alpha \in (0,a_*)$ and $L \ge 1$.
There exist a constant $K > 0$ and a period $q_* \ge 1$ such that
\[
\mathcal{A}^{(p,q)} \le
K \exp \left( -\frac{\pi \alpha}{|p/q-1/2|} \right),
\]
for all relatively prime integers $p$ and $q$ such that
$1 \le |2p-q|\le L$ and $q\ge q_*$.
\end{thm}

\section{Proof of Theorem~\ref{thm:Remainder}}
\label{sec:ProofRemainder}

\subsection{Spaces, norms, and projections}

Let $\mathcal{X}_{a,b}$, with $a>0$ and $b>0$,
be the space of all analytic functions $g$ defined on the open set
\[
D_{a,b} =
\left\{ (x,y) \in (\Cset / \Zset) \times \Cset :
|\Im x| < a, \; |y| < b \right\}
\]
with bounded Fourier norm
\[
\Fnorm{g}_{a,b} =
\sum_{k\in\Zset} \Snorm{\hat{g}_k}_b \rme^{2\pi |k| a},
\]
where $\hat{g}_k(y)$ denotes the $k$-th Fourier coefficient of the
$1$-periodic function $g(\cdot,y)$ and
\[
\Snorm{\hat{g}_k}_b = \sup \left\{ |\hat{g}_k(y)| : y \in B_b \right\}
\]
denotes its sup-norm on the complex open ball $B_b = \{ y \in \Cset : |y| < b \}$.
Let
\[
\Snorm{g}_{a,b} =
\sup \left\{ |g(x,y)| : (x,y) \in D_{a,b} \right\}
\]
be the sup-norm on $D_{a,b}$.

Let $\mathcal{X}_{a,b,m}$ be the space of all vectorial functions
$G : D_{a,b} \to \Cset^2$
of the form
\[
G(x,y)=(y^m g_1(x,y), y^{m+1}g_2(x,y))
\]
such that $g_1,g_2 \in \mathcal{X}_{a,b}$.
The space $\mathcal{\mathcal{X}}_{a,b,m}$ is a Banach space with the Fourier norm
$\Fnorm{G}_{a,b,m} =
\max \left\{\Fnorm{g_1}_{a,b}, \Fnorm{g_2}_{a,b}\right\}.$
The sup-norm $\Snorm{\cdot}_{a,b,m}$ on $\mathcal{\mathcal{X}}_{a,b,m}$
is defined analogously.

Let $g^*_2(y) = \int_0^1 g_2(x,y) \rmd x$ be the average of $g_2(x,y)$.
Let $\mathcal{X}_{a,b,m} = \mathcal{X}_{a,b,m}^* \oplus \mathcal{X}_{a,b,m}^\bullet$ be the
direct decomposition where $\mathcal{X}_{a,b,m}^*$ is the vectorial subspace
of the elements of the form $G^*(x,y) = (0,y^{m+1} g_2^*(y))$,
whereas $\mathcal{X}_{a,b,m}^\bullet$ is the one of the elements with $g^*_2(y) = 0$.
Let $\pi^*: \mathcal{X}_{a,b,m} \to \mathcal{X}_{a,b,m}^*$ and
$\pi^\bullet: \mathcal{X}_{a,b,m} \to \mathcal{X}_{a,b,m}^\bullet$ be the
associated projections.
Thus, any $G \in \mathcal{X}_{a,b,m}$ can be decomposed as
$G = G^* + G^\bullet$, where
\[
G^* = \pi^*(G) = (0,y^{m+1} g^*_2(y)) \in \mathcal{X}_{a,b,m}^*,\qquad
G^\bullet = \pi^\bullet(G) \in \mathcal{X}_{a,b,m}^\bullet.
\]
Obviously, $\Fnorm{G^*}_{a,b,m} \le \Fnorm{G}_{a,b,m}$ and
$\Fnorm{G^\bullet}_{a,b,m} \le \Fnorm{G}_{a,b,m}$

We will always denote the scalar functions in $\mathcal{X}_{a,b}$ with lower-case letters,
and the vectorial functions in $\mathcal{X}_{a,b,m}$ with upper-case letters.
Asterisk and bullet superscripts in upper-case letters stand for the
$\pi^*$-projections and $\pi^\bullet$-projections of vectorial functions
in $\mathcal{X}_{a,b,m}$, respectively.
Asterisk superscripts in lower-case letters
denote averages of scalar functions in $\mathcal{X}_{a,b}$.
We will always write the couple of scalar functions associated to
any given vectorial function of $\mathcal{X}_{a,b,m}$ with the corresponding
lower-case letter and the subscripts $j=1,2$.
Hat symbols denote Fourier coefficients.

\subsection{The averaging and the iterative lemmas}

Henceforth, let $A(x,y) = (x+y,y)$ be the integrable twist map
introduced in~(\ref{eq:IntegrableMap}).
Let $F = F_2$ be a map satisfying the properties listed in Lemma~\ref{lem:Order2},
so $F_2 = A + G_2$ for some $G_2 \in \mathcal{X}_{a_2,b_2,2}$,
where $a_2 = a_*$ is the analyticity strip width in the angular variable $x$,
and $b_2 = b_*$ is the analyticity radius in $y$.
Hence, $F_2$ is a perturbation of $A$ of order two.
The following lemma allows us to increase that order as much as we want
by simply losing as little analyticity strip width as we want.
It is based on classical averaging methods.
In particular, we see that $F$ is a perturbation beyond all order of $A$.

\begin{lem}[Averaging Lemma]\label{lem:BeyondAllOrders}
Let $F_2 = A + G_2$, with $G_2 \in \mathcal{X}_{a_2,b_2,2}$ and $a_2 > 0$ and $b_2 > 0$,
be a real analytic map with the intersection property on the cylinder
$\Tset \times (-b_2,b_2)$.
Let $m \ge 3$ be an integer.
Let $a_m$ be any analyticity strip width such that $a_m \in (0,a_2)$.

There exist an analyticity radius $b_m \in (0,b_2)$ and
a change of variables of the form $\Phi_m = \Identity + \Psi_m$
for some $\Psi_m \in \mathcal{X}_{a_m,b_m,1}$
such that the transformed map
$F_m = \Phi_m^{-1} \circ F_2 \circ \Phi_m$
is real analytic, has the intersection property on the cylinder
$\Tset \times (-b_m,b_m)$, and has the form
$F_m = A + G_m$ for some $G_m \in \mathcal{X}_{a_m,b_m,m}$.

Besides, the change of variables $\Phi_m$ is close to the
identity on $\Tset \times (-b_m,b_m)$. That is,
\begin{equation}\label{eq:ChangePhi_m}
\Phi_m(x,y) = \left( x+\Or(y),y+\Or(y^2) \right),\qquad
\det[\Phi_m(x,y)] = 1 + \Or(y),
\end{equation}
uniformly for all $(x,y) \in \Tset \times (-b_m,b_m)$.
\end{lem}

\proof
The change $\Phi_m$ is the composition of $m-2$ changes of the form
\[
\tilde{\Phi}_l = \Identity + \tilde{\Psi}_l,\qquad
\tilde{\Psi}_l \in \mathcal{X}_{a_l,b_l,l-1},\qquad
2 \le l < m,
\]
where $a_l = a_2 - (l-2)\epsilon$, $\epsilon = (a_2 - a_m)/(m-2)$,
$(b_l)_{2 \le l < m}$ is a positive decreasing sequence,
and $\tilde{\Psi}_l$ is constructed as follows to increase
the order of the perturbation from $l$ to $l+1$.

Let us suppose that we have a real analytic map with
the intersection property on $\Tset \times (-b_l,b_l)$ of the form
$F_l = A + G_l$, for some $G_l \in \mathcal{X}_{a_l,b_l,l}$
with $a_l,b_l > 0$ and $l \ge 2$.

We begin with a formal computation.
We write
\[
F_l(x,y)=
\left(x+y+y^l h_1(x)+\Or(y^{l+1}), y+y^{l+1} h_2(x)+\Or(y^{l+2})\right),
\]
where the functions $h_1(x)$ and $h_2(x)$ are 1-periodic and analytic
on the open complex strip $\{ x \in \Cset/\Zset: |\Im x| < a_l \}$.
We will see, by using an \emph{a posteriori} reasoning, that the intersection property
implies that $h_2(x)$ has zero average; that is,
$h_2^* = \int_0^1 h_2(x) \rmd x = 0$.
Nevertheless, we can not prove it yet.
Thus, we will keep an eye on $h_2^*$ in what follows.

If we take the change of variables
\begin{equation}\label{eq:PhiBeyondAllOrders}
\tilde{\Phi}_l(x,y) = \left(x + y^{l-1} \psi_1(x), y + y^l \psi_2(x) \right)
\end{equation}
for some functions $\psi_1(x)$ and $\psi_2(x)$,
then, after a straightforward computation, the map
$F_{l+1} = (\tilde{\Phi}_l)^{-1} \circ F_l \circ \tilde{\Phi}_l$
has the form
\[
F_{l+1}(x,y)=
\left(x+y+ y^l k_1(x)+\Or(y^{l+1}),
y+y^{l+1} k_2(x)+\Or(y^{l+2})\right),
\]
with $k_1=\psi_2+h_1-\psi_1'$ and $k_2=h_2-\psi_2'$.
Therefore, we take
\[
\psi_2(x) = \int_0^x \big(h_2(s) - h_2^* \big) \rmd s - h_1^*, \qquad
\psi_1(x) = \int_0^x \big( \psi_2(s)+h_1(s)\big) \rmd s,
\]
so that $k_1(x) = 0$ and $k_2(x) = h_2^*$.
These functions $\psi_2(x)$ and $\psi_1(x)$ are 1-periodic,
because $h_2(x)-h_2^*$ and $\psi_2(x)+h_1(x)$ have zero average.
Besides, $\psi_1(x)$ and $\psi_2(x)$ are analytic in the
open complex strip $\{ x\in \Cset/\Zset : |\Im x| < a_l \}$.
Indeed, $\tilde{\Phi}_l = \Identity + \tilde{\Psi}_l$
with $\tilde{\Psi}_l \in \mathcal{X}_{a_l,b_l,l-1}$.

Next, we control the domain of definition of the map $F_{l+1}$.
The inverse change is
\[
(\tilde{\Phi}_l)^{-1}(x,y)=
\left( x - y^{l-1} \psi_1(x) + \Or(y^l),
       y - y^l \psi_2(x) + \Or(y^{l+1}) \right).
\]
Thus, the maps $\tilde{\Phi}_l$, $F_l$, and $(\tilde{\Phi}_l)^{-1}$
have the form $(x,y)\mapsto(x+\Or(y),y+\Or(y^2))$,
since $l \ge 2$.
Consequently, if $b_{l+1} \le b_l/2$ is small enough, then
\[
D_{a_{l+1},b_{l+1}} \stackrel{\tilde{\Phi}_l}{\longrightarrow}
D_{a_l-2\epsilon/3, 4 b_{l+1}/3} \stackrel{F_l}{\longrightarrow}
D_{a_l -\epsilon/3,5 b_{l+1}/3} \stackrel{(\tilde{\Phi}_l)^{-1}}{\longrightarrow}
D_{a_l,2 b_{l+1}} \subset D_{a_l,b_l},
\]
so $F_{l+1} = (\Phi_l)^{-1} \circ F_l \circ \Phi_l$
is well-defined on $D_{a_{l+1}, b_{l+1}}$.
Now, let us check that $h_2^* = 0$.
At this moment, we only know that
\[
F_{l+1}(x,y)=
\left( x + y + \Or(y^{l+1}), y + y^{l+1} h_2^* + \Or(y^{l+2})\right),
\]
since the change of variables has not eliminated the average $h_2^*$.
The map $F_{l+1}$ has the intersection property on the cylinder
$\Tset \times (-b_{l+1},b_{l+1})$,
because the intersection property is preserved by changes of variables.
This implies that $h_2^* = 0$.
On the contrary,
the image of the loop $\Tset \times \{ y_0 \}$ does not intersect itself
when $0 < y_0 \ll 1$.

Finally, properties~(\ref{eq:ChangePhi_m}) follow directly from the fact
that we have performed a finite number of changes, all of them satisfying
these same properties.
\qed

Next, the following theorem provides the exponentially small bound for the
$\pi^\bullet$-projection of the residue provided an initial order $m$ big enough.
It is the main tool to prove Theorem~\ref{thm:Remainder}.

\begin{thm}\label{thm:Technical}
Let $m \ge 6$ be an integer, $\bar{a} > 0$, $\bar{d} > 0$, and $r \in (0,1)$.
There exist constants $\bar{b} = \bar{b}(m,\bar{a},\bar{d},r) > 0$
and $c_j = c_j(r) > 0$, $j=1,2,3$, such that, if
\begin{equation}\label{eq:barFv}
\bar{F} = A + \bar{G},\quad
\bar{G} \in \mathcal{X}_{\bar{a},\bar{b},m},\quad
\bar{d}^* = \Fnorm{\pi^*(\bar{G})}_{\bar{a},\bar{b},m},\quad
\bar{d}^\bullet = \Fnorm{\pi^\bullet(\bar{G})}_{\bar{a},\bar{b},m},
\end{equation}
and
\[
0 < \breve{a} < \bar{a},\qquad
0 < \breve{b} \le \bar{b}\sqrt{r},\qquad
\bar{d}^* + (1 + c_2) \bar{d}^\bullet \le \bar{d},
\]
then there exists a change of variables
$\breve{\Phi} = \Identity + \breve{\Psi}$ satisfying the following properties:
\begin{enumerate}
\item
$\breve{\Psi} \in \mathcal{X}_{\breve{a},\breve{b},m-1}$ with
$\Snorm{\breve{\Psi}}_{\breve{a},\breve{b},m-1} \le c_1 \bar{d}^\bullet$; and
\item
The transformed map
$\breve{F} = \breve{\Phi}^{-1} \circ \bar{F} \circ \breve{\Phi}$
is real analytic, has the intersection property on the cylinder
$\Tset \times (-\breve{b},\breve{b})$,
and has the form $\breve{F} = A + \breve{G}$,
$\breve{G} \in \mathcal{X}_{\breve{a},\breve{b},m}$,
\[
\Fnorm{\pi^*(\breve{G})}_{\breve{a},\breve{b},m} \le
\bar{d}^* + c_2 \bar{d}^\bullet,
\qquad
\Fnorm{\pi^\bullet(\breve{G})}_{\breve{a},\breve{b},m} \le
c_3 \rme^{-2\pi r (\bar{a}-\breve{a})/\breve{b}} \bar{d}^{\bullet}.
\]
\end{enumerate}
\end{thm}

Theorem~\ref{thm:Technical} is proved in Subsection~\ref{sec:ProofThmTechnical}.

In order to present the main ideas of the proof,
let us try to completely get rid of the remainder of the map of the form
$F = A + G$, for some $G \in \mathcal{X}_{a,b,m}$, with a change of variables
of the form $\Phi = \Identity + \Psi$,
for some $\Psi \in \mathcal{X}_{a,b,m-1}$.
Concretely, we look for $\Phi$ such that
$A = \Phi^{-1} \circ F \circ \Phi$,
or, equivalently, we look for $\Psi$ such that
\[
\Psi \circ A - A \Psi= G \circ (\Identity+\Psi).
\]
It is not possible to solve this equation in general.
Instead, we consider the linear equation
\[
\Psi \circ A - A \Psi = G.
\]
This vectorial equation reads as
\[
\left\{\begin{array}{l}
\psi_1(x+y,y)-\psi_1(x,y)  =  y \left(\psi_2(x,y) + g_1(x, y)\right), \\
\psi_2(x+y,y)-\psi_2(x,y)  =  y g_2(x,y).
\end{array}\right.
\]
Therefore, we need to solve two linear equations of the form
\begin{equation}\label{eq:LinearSingle}
\psi(x+y,y)-\psi(x,y) = y g(x,y),
\end{equation}
where $g \in \mathcal{X}_{a,b}$ is known.
If the average of $g(x,y)$ is different from zero:
$g^*(y) = \hat{g}_0(y) \neq 0$,
then this equation can not be solved.
Besides, it is a straightforward computation to check that, if $\hat{g}_0(y) = 0$,
the formal solution of this equation in the Fourier basis is
\begin{equation}\label{eq:FourierCoefficient}
\hat{\psi}_k(y) = \frac{y}{\rme^{2\pi k y \rmi}-1} \hat{g}_k(y),
\qquad \forall k \neq 0,
\end{equation}
whereas the zero-th coefficient $\hat{\psi}_0(y)$ can be chosen arbitrarily.
From~(\ref{eq:FourierCoefficient}), it is clear that~(\ref{eq:LinearSingle})
can not be solved unless $g$ has only a finite number of
harmonics and zero average.
For this reason, given a function $g(x,y)$ with zero average,
we define its $K$\emph{-cut off} as
\begin{equation}\label{eq:KCutOff}
g^{<K}(x,y) = \sum_{|k| < K} \hat{g}_k(y)\rme^{2\pi kx\rmi}.
\end{equation}
Let $K$ be such that $|2\pi k y| < 2\pi$ for all $|y| < b$ and $|k| < K$.
Hence, we will take $K= s/b$ for some fixed $s\in(0,1)$,
and we will actually solve truncated linear equations of the form
\begin{equation}\label{eq:TruncatedLinearSingle}
\psi(x+y,y)-\psi(x,y) = y g^{< K}(x,y).
\end{equation}
The Fourier norm is specially suited to analyse
this kind of equations; see Lemma~\ref{lem:LinearSingle}.

Summarizing these ideas,
we look for a change of variables of the form $\Phi=\Identity+\Psi$,
where $\Psi$ satisfies the truncated linear vectorial equation
\begin{equation}\label{eq:TruncatedLinearDouble}
\Psi \circ A - A \Psi = (G^\bullet)^{<K},
\end{equation}
where $(G^\bullet)^{<K}$ denotes the $K$-cut off of
$G^\bullet = \pi^\bullet(G)$.
The average of the first component of $(G^\bullet)^{<K}$ may be non-zero.
Equation~(\ref{eq:TruncatedLinearDouble}) is studied in
Lemma~\ref{lem:LinearDouble}.
This close to the identity change of variables $\Phi=\Identity+\Psi$
does not completely eliminate the remainder.
However, if $b$ is small enough,
it reduces the size of the $\pi^\bullet$-projection of the remainder
as the following lemma shows.

\begin{lem}[Iterative lemma]\label{lem:Iterative}
Let $m \ge 6$ be an integer, $\bar{a} > 0$, $\bar{d} > 0$, $\mu > 0$, and $\rho \in (0,1)$.
There exist a constant $\bar{b} > 0$ such that if
\[
F = A + G, \quad
G\in\mathcal{X}_{a,b,m},\quad
d^* = \Fnorm{\pi^*(G)}_{a,b,m},\quad
d^\bullet = \Fnorm{\pi^\bullet(G)}_{a,b,m},
\]
with
\[
0 < a \le \bar{a}, \qquad
0<b \le \min\{\bar{b},a/6\},\qquad
d^* + d^\bullet \le \bar{d},
\]
then there exists a change of variables $\Phi = \Identity+\Psi$
satisfying the following properties:
\begin{enumerate}
\item
$\Psi \in \mathcal{X}_{a,b,m-1}$ is a solution of the truncated linear
equation~(\ref{eq:TruncatedLinearDouble}) such that
\[
\Snorm{\Psi}_{a,b,m-1} \le \Fnorm{\Psi}_{a,b,m-1} \le \Omega d^\bullet,
\]
where $\Omega = \Omega(\sqrt{\rho})$ is defined in Lemma~\ref{lem:LinearDouble}; and
\item
The transformed map $\tilde{F}= \Phi^{-1}\circ F \circ \Phi$
is real analytic, has the intersection property on the cylinder
$\Tset \times (-\tilde{b},\tilde{b})$, and has the form
$\tilde{F} = A+\tilde{G}$, $\tilde{G}\in\mathcal{X}_{\tilde{a},\tilde{b},m}$,
\[
\Fnorm{\pi^*(\tilde{G})}_{\tilde{a},\tilde{b},m} \le
d^* + \rme^{-12\pi\rho} d^\bullet,\qquad
\Fnorm{\pi^\bullet(\tilde{G})}_{\tilde{a},\tilde{b},m} \le
\rme^{-12\pi \rho} d^\bullet,
\]
where $\tilde{a} = a-6b$ and $\tilde{b} = b - \mu b^2$.
\end{enumerate}
\end{lem}

\begin{remark}
If $\tilde{a} = a-6b$, then
$\rme^{-2\pi \rho (a-\tilde{a})/b} = \rme^{-12 \pi \rho}$.
\end{remark}

The proof of this lemma is found in Subsection~\ref{ssec:ProofOfIterativeLemma}.
Some technicalities in the proof require the use of the sup-norm,
which forces us to deal with both the Fourier norm and the sup-norm.
The relations between them are stated in Lemma~\ref{lem:NormBounds}.

Finally, Theorem~\ref{thm:Technical} is obtained by means of a finite sequence of
changes of variables like the ones described in the iterative lemma.
We want to perform as many of such changes as possible because each change reduces the
size of the $\pi^\bullet$-projection of the remainder by the factor $\rme^{-12\pi \rho}$.
Since the loss of analyticity in the angular variable is $6b = \Or(b)$,
then we can at most perform a number $\Or(1/b)$ of such changes.
This idea goes back to Neishtadt~\cite{Neishtadt1984}.

The intersection property is used neither in the proof of the iterative lemma
nor in the proof of Theorem~\ref{thm:Technical},
but will be essential to control the size of the $\pi^*$-projections
of the remainders in terms of the size of their $\pi^\bullet$-projections
later on.

\subsection{Technical lemmas}

\begin{lem}\label{lem:NormBounds}
Let $0 < \alpha < \min\{a,1/2\pi\}$, $b > 0$, and $g \in \mathcal{X}_{a,b}$.
Let $g^{\ge K}=g-g^{<K}$, with $g^{<K}$ the $K$-cut off of $g$,
defined in~(\ref{eq:KCutOff}).
Then:
\begin{enumerate}
\item
$\Fnorm{g^{< K}}_{a,b} \leq \Fnorm{g}_{a,b}$,
\item
$\Fnorm{g^{\ge K}}_{a-\alpha,b} \leq \rme^{-2\pi K\alpha} \Fnorm{g}_{a,b}$,
\item
$\Snorm{g}_{a,b}\leq\Fnorm{g}_{a,b}$, and
\item
$\Fnorm{g}_{a-\alpha,b} \leq {\alpha}^{-1}\Snorm{g}_{a,b}$.
\end{enumerate}
If $m \in \Nset$, then these bounds also hold for any vectorial
function $G \in \mathcal{X}_{a,b,m}$.
\end{lem}

\proof
First, the Fourier norm of $g^{<K}$ is a partial sum of the Fourier norm of $g$.
Second,
$\Fnorm{g^{\ge K}}_{a-\alpha,b} =
\sum_{|k|\ge K} \Snorm{\hat{g}_k}_b \rme^{2\pi |k| (a-\alpha)} =
\rme^{-2\pi K \alpha} \sum_{|k|\ge K} \Snorm{\hat{g}_k}_b \rme^{2\pi |k| a} \le
\rme^{-2\pi K \alpha} \Fnorm{g}_{a,b}$.
Third,
$|g(x,y)| \le
 \sum_{k\in \Zset} |\hat{g}_k(y)| |\rme^{2\pi k x \rmi}| \le
 \sum_{k\in \Zset} \Snorm{\hat{g}_k}_b  \rme^{2\pi |k|a} =
 \Fnorm{g}_{a,b}$,
for all $(x,y) \in D_{a,b}$.
Fourth, we recall that the Fourier coefficients of the analytic function $g$
satisfy the inequality $\Snorm{\hat{g}_k}_b \le \rme^{-2\pi|k|a} \Snorm{g}_{a,b}$
for all $k \in \Zset$.
Hence,
\[
\Fnorm{g}_{a-\alpha,b} =
\sum_{k\in\Zset} \Snorm{\hat{g}_k}_b \rme^{2\pi |k| (a-\alpha)} \le
2 \Snorm{g}_{a,b} \sum_{k\ge 0} \rme^{-2\pi k\alpha} \le
\alpha^{-1}\Snorm{g}_{a,b},
\]
where we have used that
$\sum_{k \ge 0} \rme^{- kt} = (1-\rme^{-t})^{-1}\le \rme/t < \pi/t$ for all $t\in(0,1)$.
The last part follows from the definition of the norms $\Fnorm{\cdot}_{a,b,m}$
and $\Snorm{\cdot}_{a,b,m}$.
\qed

\begin{lem}\label{lem:LinearSingle}
If $s \in (0,1)$, $K = s/b$, and $g \in \mathcal{X}_{a,b}$ is
a function with zero average,
then the truncated linear equation~(\ref{eq:TruncatedLinearSingle})
has a unique solution $\psi \in \mathcal{X}_{a,b}$
with zero average and $\Fnorm{\psi}_{a,b} \le \omega \Fnorm{g}_{a,b}$,
where
\begin{equation}\label{eq:omega}
\omega = \omega(s) =
\frac{1}{2\pi} \cdot \max_{|z| \le 2\pi s} \left|\frac{z}{\rme^z-1} \right|.
\end{equation}
\end{lem}

\proof
The Fourier coefficients of $\psi$ must satisfy~(\ref{eq:FourierCoefficient}).
We note that $\omega < \infty$ for all $s \in (0,1)$,
since the function $z/(\rme^z-1)$ is analytic on the open ball $|z| < 2 \pi$.
Moreover,
\[
\Snorm{\hat{\psi}_k}_b \le
\left( \max_{|y| \le b} \left|{\frac{y}{\rme^{2\pi k y \rmi}-1}}\right|\right)
\Snorm{\hat{g}_k}_b \le
 \frac{\omega}{|k|} \Snorm{\hat{g}_k}_b \le
 \omega \Snorm{\hat{g}_k}_b,
\]
for all $0 < |k| < K=s/b$.
Finally, we recall that $\hat{\psi}_0(y) \equiv 0$. Then we obtain that
$\Fnorm{\psi}_{a,b} =
\sum_{0 < |k| < K} \Snorm{\hat{\psi}_k}_b \rme^{2\pi |k| a} \le
\omega \sum_{k \in \Zset} \Snorm{\hat{g}_k}_b \rme^{2\pi |k| a} =
\omega \Fnorm{g}_{a,b}$.
\qed

\begin{remark}
We will denote by  $\psi = \mathcal{G}_K(g^{<K})$
the linear operator that
sends the independent term $g^{<K}$ of the truncated linear
equation~(\ref{eq:TruncatedLinearSingle}) to the solution $\psi$
with zero average.
Note that the solution $\psi$ has no harmonics of order $\ge K$.
\end{remark}

\begin{lem}\label{lem:LinearDouble}
If $m \ge 1$, $s \in (0,1)$, $K = s/b$, and $G\in\mathcal{X}_{a,b,m}$,
then the truncated linear equation~(\ref{eq:TruncatedLinearDouble})
has a solution $\Psi\in\mathcal{X}_{a,b,m-1}$
such that
\[
\Fnorm{\Psi}_{a,b,m-1}\le \Omega \Fnorm{G^\bullet}_{a,b,m},
\]
where $G^\bullet = \pi^{\bullet}(G)$,
$\Omega = \Omega(s) = ( \omega(s) + 1 ) \max\{1,\omega(s)\}$,
and $\omega(s)$ is defined in~(\ref{eq:omega}).
\end{lem}

\proof
Let $G = (y^m g_1, y^{m+1}g_2)$ and $\Psi = (y^{m-1} \psi_1, y^m\psi_2)$.
Then the vectorial equation $\Psi \circ A - A \Psi = (G^\bullet)^{<K}$ reads as
\[
\left\{\begin{array}{l}
\psi_1(x+y,y)-\psi_1(x,y) =  y \left( \psi_2(x,y) + g_1^{<K}(x,y) \right),\\
\psi_2(x+y,y)-\psi_2(x,y) = y (g_2^{<K} (x,y) - g_2^*(y) ).
\end{array}\right.
\]
Let $\psi_2=\mathcal{G}_K(g_2^{<K}-g_2^*)-g_1^*$ and
$\psi_1=\mathcal{G}_K(\psi_2+g_1^{<K})$.
These operations are well-defined since
both $g_2-g_2^*$ and $\psi_2+g_1$ have zero average.
As for the bounds,
\begin{eqnarray*}
\Fnorm{\psi_2}_{a,b} & \le &
\omega \Fnorm{g_2 - g_2^*}_{a,b} + \Fnorm{g_1}_{a,b}
\le \Omega \Fnorm{G^\bullet}_{a,b,m},\\
\Fnorm{\psi_1}_{a,b} & \le &
\omega \Fnorm{\psi_2 + g_1^{<K}}_{a,b} \\
& \le &
\omega \Fnorm{\psi_2 + g_1^*}_{a,b} +
\omega \Fnorm{g_1^{<K} - g_1^*}_{a,b} \\
& \le &
\omega^2 \Fnorm{g_2 - g_2^*}_{a,b} + \omega \Fnorm{g_1}_{a,b} \le
\Omega \Fnorm{G^\bullet}_{a,b,m},
\end{eqnarray*}
where we have used Lemma~\ref{lem:LinearSingle}.
\qed

\begin{lem}\label{lem:Difference}
Let $l,n\in\Nset$, $0<\alpha<\min\{a/3,1/2\pi\}$, $0<\beta<b/2$,
$c_1,c_2>0$, and $c = c_1+c_2$, such that
\begin{eqnarray}\label{eq:abcd}
b + b^n c <  \alpha + \beta,\qquad  b^{n+1} c < \beta.
\end{eqnarray}
Let $M = M(\alpha, b, \beta, c_1, c_2, l, n) =
(1 + b^n c )^{l+1} b^{l-1}
\left( \alpha^{-1} + b\beta^{-1} + l+1 \right)$.
Let $\Delta\in\mathcal{X}_{a,b,l}$.
Let $\Gamma_1, \Gamma_2 \in \mathcal{X}_{a-2\alpha, b-2\beta,n}$ with
$\Fnorm{\Gamma_j}_{a-2\alpha, b-2\beta,n} \le c_j$.
Let $L(x,y) = (x + \eta y,y)$ with $|\eta| \le 1$.
Then,
\begin{enumerate}
\item
$\Delta\circ(L+\Gamma_1) - \Delta\circ(L+\Gamma_2) \in \mathcal{X}_{a-3\alpha,b-2\beta,n+l}$,
\item
$\Snorm{\Delta\circ(L+\Gamma_1) - \Delta\circ(L+\Gamma_2)}_{a-2\alpha,b-2\beta,n+1} \le
M \Snorm{\Delta}_{a,b,l} \Snorm{\Gamma_1 - \Gamma_2}_{a-2\alpha,b-2\beta,n}$, and
\item
$\Fnorm{\Delta\circ(L+\Gamma_1) - \Delta\circ(L+\Gamma_2)}_{a-3\alpha,b-2\beta,n+1} \le
M \alpha^{-1} \Fnorm{\Delta}_{a,b,l} \Fnorm{\Gamma_1 - \Gamma_2}_{a-2\alpha,b-2\beta,n}$.
\end{enumerate}
\end{lem}

\proof
Let $\Gamma = \Gamma_1 - \Gamma_2$.
Then $\Snorm{\Gamma}_{a-2\alpha,b-2\beta,n} \le
 \Fnorm{\Gamma}_{a-2\alpha,b-2\beta,n} \le c$
and
\[
\Delta\circ(L+\Gamma_1) - \Delta\circ(L+\Gamma_2) =
\int_0^1 \big( \rmD\Delta\circ(L+t\Gamma) \big) \cdot \Gamma\ \rmd t.
\]
Let $(x_t,y_t) =
(L+t\Gamma)(x,y) =
(x+\eta y+ t y^n \gamma_1(x,y), y+ t y^{n+1} \gamma_2(x,y))$,
with $t \in [0,1]$ fixed.
We deduce from conditions~(\ref{eq:abcd}) that $(x_0,y_0)$ and $(x_1,y_1)$ belong to
$D_{a-\alpha,b-\beta}$ for all $(x,y)\in D_{a-2\alpha, b-2\beta}$.
Therefore, $(x_t,y_t)\in D_{a-\alpha,b-\beta}$ by convexity of the domain,
and so, the composition $\Delta\circ(L+t\Gamma)$ is well-defined on the domain
$D_{a-2\alpha, b-2\beta}$.

A simple computation shows that the product
$\rmD \Delta(x_t,y_t)\cdot \Gamma(x,y)$ is equal to
\begin{eqnarray*}
\fl
\left(
\begin{array}{c}
y_t^{l-1}y^{n+1} \left( y^{-1}y_t\partial_1 \delta_1 (x_t,y_t)\gamma_1(x,y) +
\left( y_t\partial_2 \delta_1(x_t,y_t) + l \delta_1(x_t,y_t) \right) \gamma_2(x,y)\right)\\
y_t^{l}y^{n+1} \left( y^{-1}y_t\partial_1 \delta_2(x_t,y_t) \gamma_1(x,y) +
\left(y_t\partial_2 \delta_2(x_t,y_t)  +(l+1)\delta_2(x_t,y_t)  \right)\gamma_2(x,y)\right)
\end{array}
\right).
\end{eqnarray*}
Let us bound the elements above.
On the one hand,
$|\gamma_i(x,y)|\le \Snorm{\Gamma}_{a-2\alpha,b-2\beta,n} \le c$,
$|y_t|= |y + t y^{n+1}\gamma_2(x,y)|\le ( 1 + b^n c) |y|$, and
$|y|<b$ for all $(x,y)\in D_{a-2\alpha, b-2\beta}$.
On the other hand, $|\delta_i(x_s,y_s)|\le \Snorm{\Delta}_{a,b,l}$
and the Cauchy estimates imply that
\begin{eqnarray*}
|\partial_1 \delta_i(x_t,y_t)| \le \alpha^{-1} \Snorm{\Delta}_{a,b,l}, \qquad
|\partial_2 \delta_i(x_t,y_t)| \le \beta^{-1} \Snorm{\Delta}_{a,b,l}.
\end{eqnarray*}
From the previous bounds and
the definitions of both norms,
we deduce that
\[
\Snorm{(\rmD\Delta\circ(L+t\Gamma) ) \cdot \Gamma }_{a-2\alpha,b-2\beta,n+l} \le
M'\Snorm{\Delta}_{a,b,l} \Snorm{\Gamma}_{a-2\alpha,b-2\beta,n},
\]
for all $t\in[0,1]$, where $M'=(1+b^nc)^{l+1} \left( \alpha^{-1} + b\beta^{-1} + l+1 \right)$.
Thus,
\begin{eqnarray*}
\fl
\Snorm{\Delta\circ(L+\Gamma_1) - \Delta\circ(L+\Gamma_2)}_{a-2\alpha,b-2\beta,n+l}
&\le&
M'\Snorm{\Delta}_{a,b,l} \Snorm{\Gamma}_{a-2\alpha,b-2\beta,n}, \\
\fl
\Fnorm{\Delta\circ(L+\Gamma_1) - \Delta\circ(L+\Gamma_2)}_{a-3\alpha,b-2\beta,n+l}
&\le&
\alpha^{-1} M'\Fnorm{\Delta}_{a,b,l} \Fnorm{\Gamma}_{a-2\alpha,b-2\beta,n}.
\end{eqnarray*}
This proves the first item.
The other items follow from the bounds
$\Snorm{\cdot}_{a,b,n+1} \le b^{l-1}\Snorm{\cdot}_{a,b,n+l}$ and
$\Fnorm{\cdot}_{a,b,n+1} \le b^{l-1}\Fnorm{\cdot}_{a,b,n+l}$,
since $M=M'b^{l-1}$.
\qed

\begin{lem}\label{lem:Inverse}
Let $n \in \Nset$, $0<\alpha<\min\{a/3,1/2\pi\}$,
$0<\beta<b/2$, and $p>0$ such that  conditions~(\ref{eq:abcd}) hold with $c=2p$.
Let $\Phi=\Identity+\Psi$, with $\Psi \in \mathcal{X}_{a,b,n}$
and $\Fnorm{\Psi}_{a,b,n} \le p$.
Then $\Phi(D_{a',b'})\subset D_{a'+\alpha,b'+\beta}$
for all $0 < a' \le a$ and  $0 < b' \le b$.

Let $M_* = M(\alpha,b,\beta,p,p,n,n)$, where $M$ is
defined in Lemma~\ref{lem:Difference}.
If $M_*bp < 1$,
then $\Phi$ is invertible and the inverse change $\Phi^{-1}$ satisfies
the following properties:
\begin{enumerate}
\item
$\Phi^{-1}=\Identity+\Upsilon$ for some
$\Upsilon \in \mathcal{X}_{a-2\alpha, b-2\beta,n}$ such that
$\Snorm{\Upsilon}_{a-2\alpha,b-2\beta,n} \le \Snorm{\Psi}_{a,b,n}$,
\item
$\Phi^{-1}(D_{a',b'})\subset D_{a'+\alpha,b'+\beta}$
for all $0 < a' \le a - 2\alpha$ and $0 < b' \le b-2\beta$, and
\item
$\Fnorm{\Upsilon + \Psi}_{a-3\alpha,b-2\beta,n+1} \le
M_*\alpha^{-1} \Fnorm{\Psi}^2_{a,b,n}$.
\end{enumerate}
\end{lem}

\proof
Note that $\Snorm{\Psi}_{a,b,n} \le \Fnorm{\Psi}_{a,b,n} \le p$.
Conditions~(\ref{eq:abcd}) imply that
$b^n p < \alpha$ and $b^{n+1} p < \beta$.
Therefore, $\Phi(D_{a',b'})\subset D_{a'+\alpha,b'+\beta}$
for all $0 < a' \le a$ and $0 < b' \le b$.
Analogously, if $\Upsilon\in\mathcal{X}_{a-2\alpha, b-2\beta, n}$ and
$\Snorm{\Upsilon}_{a-2\alpha, b-2\beta, n} \le p$,
then $(\Identity+\Upsilon)(D_{a',b'})\subset D_{a'+\alpha,b'+\beta}$
for all $0 < a' \le a - 2\alpha$, $0 < b' \le b-2\beta$.

We denote by $\mathcal{B}$ the closed ball in
$\mathcal{X}_{a-2\alpha,b-2\beta,n}$ of radius $p$ in the sup-norm.
Let us prove that the functional $\mathcal{P}:\mathcal{B} \to \mathcal{B}$,
$\mathcal{P}(\Upsilon) = - \Psi \circ (\Identity+\Upsilon)$,
is a well-defined contraction with Lipschitz constant
\begin{equation}\label{eq:LipschitzConstant}
\Lipschitz \mathcal{P} \le
 M_* b  \Snorm{\Psi}_{a,b,n} \le M_* b p < 1.
\end{equation}
First, we observe that
\begin{equation}\label{eq:Pnorm}
\Snorm{\mathcal{P}(\Upsilon)}_{a-2\alpha,b-2\beta,n} \le \Snorm{\Psi}_{a,b,n} \le p,
\qquad \forall \Upsilon\in\mathcal{B},
\end{equation}
so $\mathcal{P} (\mathcal{B}) \subset \mathcal{B}$.
Second, we bound
$\mathcal{P}(\Upsilon)-\mathcal{P}(\Xi)=
\Psi \circ(\Identity+\Xi) - \Psi \circ (\Identity + \Upsilon)$ as follows:
\begin{eqnarray}\label{eq:contraction}
\Snorm{\mathcal{P}(\Upsilon) - \mathcal{P}(\Xi)}_{a-2\alpha,b-2\beta,n}
& \le &
b \Snorm{ \mathcal{P}(\Upsilon)-\mathcal{P}(\Xi)}_{a-2\alpha, b-2\beta, n+1}
\nonumber \\
& \le &
M_* b \Snorm{\Psi}_{a, b, n} \Snorm{\Xi - \Upsilon}_{a-2\alpha, b-2\beta, n}.
\end{eqnarray}
The first inequality is direct,
and the second comes from Lemma~\ref{lem:Difference} with
$\Delta = \Psi$, $\Gamma_1 = \Xi$, $\Gamma_2 = \Upsilon$,
$L = \Identity$, $c_j=p$, and $l=n$.
This proves that $\mathcal{P}$ is a contraction
with Lipschitz constant~(\ref{eq:LipschitzConstant}).
Thus, $\mathcal{P}$ has a unique fixed point $\Upsilon \in \mathcal{B}$
which satisfies that
\[
(\Identity+\Psi) \circ (\Identity+\Upsilon) =
\Identity + \Upsilon + \Psi \circ (\Identity+\Upsilon) =
\Identity + \Upsilon - \mathcal{P}(\Upsilon) =
\Identity
 \]
on $D_{a-2\alpha,b-2\beta}$.
Therefore, the inverse map $\Phi^{-1}$ exists and equals $\Identity+\Upsilon$.
Furthermore, $\Snorm{\Upsilon}_{a-2\alpha,b-2\beta,n}\le\Snorm{\Psi}_{a,b,n}$
follows from~(\ref{eq:Pnorm}).
Finally,
\begin{eqnarray*}
\Fnorm{\Upsilon + \Psi }_{a-3\alpha,b-2\beta,n+1}  & \le &
\alpha^{-1}\Snorm{\Upsilon + \Psi }_{a-2\alpha,b-2\beta,n+1}  \\
& \le &
\alpha^{-1} \Snorm{ \mathcal{P}(\Upsilon) - \mathcal{P}(\mathbf{0}) }_{a-2\alpha,b-2\beta,n+1}\\
& \le &
M_*\alpha^{-1}\Snorm{\Psi}_{a,b,n} \Snorm{\Upsilon}_{a-2\alpha,b-2\beta,n} \\
& \le &
M_*\alpha^{-1} \Snorm{\Psi}^2_{a,b,n}\le M_*\alpha^{-1} \Fnorm{\Psi}^2_{a,b,n}.
\end{eqnarray*}
We have used the second inequality of equation~(\ref{eq:contraction})
with $\Xi=\mathbf{0}$.
\qed

\subsection{Proof of Lemma~\ref{lem:Iterative}}\label{ssec:ProofOfIterativeLemma}

We recall that $F = A + G$ with
$G = \pi^*(G) + \pi^\bullet(G) = G^* + G^\bullet \in \mathcal{X}_{a,b,m}$,
$d^* = \Fnorm{G^*}_{a,b,m}$, $d^\bullet = \Fnorm{G^\bullet}_{a,b,m}$,
$d^* + d^\bullet \le \bar{d}$, and $m \ge 6$.
Let $s = \sqrt{\rho} \in(\rho,1)$.
Let $\Omega = \Omega(s) = \Omega(\sqrt{\rho}) > 0$ be the constant introduced in Lemma~\ref{lem:LinearDouble}
and $\sigma = 1 + 2\Omega$.
Let $\Phi= \Identity + \Psi$ be the change of variables
where $\Psi\in \mathcal{X}_{a,b,m-1}$ is the solution given in Lemma~\ref{lem:LinearDouble}
of the truncated linear equation $\Psi \circ A - A\Psi = (G^\bullet)^{<K}$ with $K=s/b$,
so that
\begin{equation}\label{eq:PsiGbullet}
\Snorm{\Psi}_{a',b',m-1} \le
\Fnorm{\Psi}_{a',b',m-1} \le
\Omega \Fnorm{G^{\bullet}}_{a',b',m} \le
\Omega d^\bullet,
\end{equation}
for all $0<a'\le a$ and $0<b'\le b$.
Let $\Phi^{-1} = \Identity + \Upsilon$ be the inverse change studied
in Lemma~\ref{lem:Inverse}.
Let $\tilde{F} = \Phi^{-1} \circ F \circ \Phi$ be the transformed map.
Let $\tilde{G}  = \tilde{F} - A$ be the new remainder.

Henceforth, we will assume that $\alpha$, $b$, and $\beta$
are some positive constants such that
\begin{equation}\label{eq:Conditions}
\fl
b \le \alpha < \min\{a/6,1/2\pi\},\quad
0 < \beta < b/4, \quad
b+ \sigma b^{m-1}\bar{d}<\alpha+\beta,\quad
\sigma b^{m}\bar{d}<\min(\alpha,\beta).
\end{equation}

We split the proof in four steps.

\noindent
\emph{Step 1: Control of the domains.}
Note that $\tilde{F}(D_{a',b'}) \subset D_{a'+4\alpha,b'+3\beta}$
for all $0 < a' \le a - 4\alpha$ and $0 < b' \le b - 4\beta$.
Indeed,
\[
\tilde{F} : D_{a',b'} \stackrel{\Phi}{\longrightarrow}
D_{a'+\alpha,b'+\beta} \stackrel{F}{\longrightarrow}
D_{a'+3\alpha,b'+2\beta} \stackrel{\Phi^{-1}}{\longrightarrow}
D_{a'+4\alpha,b'+3\beta}.
\]
The behaviours of the changes $\Phi$ and $\Phi^{-1}$
follow directly from Lemma~\ref{lem:Inverse}, which can be applied since
conditions~(\ref{eq:Conditions}) are more restrictive than the ones required in
Lemma~\ref{lem:Inverse} when $p=\Omega\bar{d}$ and $n=m-1$.
We also need that
$a'+2\alpha \le a-2\alpha$ and $b'+2\beta \le b - 2\beta$ in order to control
the inverse $\Phi^{-1}$, which explains the restrictions on $a'$ and $b'$.

The behaviour of the map $F= A + G$ follows from the bound
\[
\Snorm{G}_{a,b,m} \le
\Fnorm{G}_{a,b,m} \le
\Fnorm{G^*}_{a,b,m} + \Fnorm{G^\bullet}_{a,b,m} =
d^* + d^\bullet \le \bar{d}
\]
and conditions
$b + b^m\bar{d} < 2\alpha$ and $b^{m+1}\bar{d} < \beta$, which are also a
consequence of (\ref{eq:Conditions}).

\noindent
\emph{Step 2: Decomposition of the new remainder.}
It turns out that $\tilde{G} = G^* + \sum_{j=1}^4 \tilde{G}_j$, where
\begin{equation}\label{eq:Gvj}
\begin{array}{ll}
\tilde{G}_1 = (G^\bullet)^{\ge K} = G^\bullet - (G^\bullet)^{<K}, \quad &
\tilde{G}_2= G \circ \Phi - G,\\
\tilde{G}_3= \Psi \circ A - \Psi \circ F \circ \Phi, \quad &
\tilde{G}_4= (\Upsilon+\Psi)\circ(F\circ\Phi).
\end{array}
\end{equation}
Indeed, $G^* + \tilde{G}_1 + \tilde{G}_2 = G\circ\Phi -(G^\bullet)^{<K}$ and
$\tilde{G}_3 + \tilde{G}_4 = \Psi \circ A + \Upsilon \circ (F \circ \Phi)$, so
\begin{eqnarray*}
G^* + \textstyle \sum_{j=1}^4 \tilde{G}_j & = &
G \circ \Phi + A \Psi + \Upsilon \circ (F\circ\Phi) \\
& = &
(F - A)\circ\Phi + A(\Phi - \Identity) +
(\Phi^{-1} - \Identity) \circ (F\circ\Phi) \\
& = &
\Phi^{-1} \circ F \circ \Phi - A = \tilde{G}.
\end{eqnarray*}
Finally, let $\tilde{G}^* = \pi^*(\tilde{G}) = G^* + \sum_{j=2}^4 \pi^*(G_j)$
and $\tilde{G}^\bullet = \pi^\bullet(\tilde{G}) = \sum_{j=1}^4 \pi^\bullet(G_j)$.

\noindent
\emph{Step 3: Bounds of the projections of the new remainder.}
Lemma~\ref{lem:NormBounds} and the bound~(\ref{eq:PsiGbullet})
will be used several times in what follows.
Below, we apply Lemma~\ref{lem:Difference} (twice) and
Lemma~\ref{lem:Inverse} (once).
The required hypotheses in each case are satisfied
due to conditions~(\ref{eq:Conditions}).

\begin{itemize}
\item
If $\tilde{a} < a$ and $\tilde{b} \le b$, then
\[
\Fnorm{\tilde{G}_1}_{\tilde{a},\tilde{b},m} =
\Fnorm{(G^\bullet)^{\ge K}}_{\tilde{a},\tilde{b},m} \le
\rme^{-2\pi K (a-\tilde{a})} \Fnorm{G^\bullet}_{a,\tilde{b},m} \le
\rme^{-2\pi K (a-\tilde{a})} d^\bullet.
\]

\item
If $\tilde{a} \le a-3\alpha$ and $\tilde{b} \le b - 2\beta$, then
\begin{eqnarray*}
\Fnorm{\tilde{G}_2}_{\tilde{a},\tilde{b},m} & \le &
M_2 \alpha^{-1}
\Fnorm{G}_{\tilde{a}+3\alpha,\tilde{b}+2\beta,m}
\Fnorm{\Psi}_{\tilde{a}+\alpha,\tilde{b},m-1} \\
& \le &
\Omega M_2 \alpha^{-1}
\Fnorm{G}_{\tilde{a}+3\alpha,\tilde{b}+2\beta,m}
\Fnorm{G^\bullet}_{\tilde{a}+\alpha,\tilde{b},m} \\
& \le &
\Omega M_2 \alpha^{-1} \bar{d} d^\bullet.
\end{eqnarray*}
The first inequality follows from Lemma~\ref{lem:Difference} with
$\Delta=G$, $L=\Identity$, $\Gamma_1=\Psi$,
$\Gamma_2=\mathbf{0}$, $l=m$, and $n=m-1$,
so that $M_2=M(\alpha,b,\beta,\Omega\bar{d},0,m,m-1)$.

\item
If $\tilde{a} \le a-2\alpha$ and $\tilde{b} \le b - \beta$, then
$\Fnorm{F \circ \Phi - A}_{\tilde{a},\tilde{b},m-1} \le
 \sigma \Fnorm{G}_{a,b,m}$.
Indeed, $F \circ \Phi - A=A\Psi+G\circ\Phi$ and
\begin{eqnarray*}
\Fnorm{A \Psi}_{\tilde{a},\tilde{b},m-1} & \le &
2\Fnorm{\Psi}_{\tilde{a},\tilde{b},m-1} \le
2\Omega\Fnorm{G^\bullet}_{\tilde{a},\tilde{b},m} \le
2\Omega d^\bullet \le 2\Omega \bar{d}, \\
\Fnorm{G \circ \Phi}_{\tilde{a},\tilde{b},m-1} & \le &
\alpha^{-1} \Snorm{G \circ \Phi}_{\tilde{a}+\alpha,\tilde{b},m-1} \le
\alpha^{-1} \Snorm{G}_{\tilde{a}+2\alpha,\tilde{b}+\beta,m-1} \\
& \le &
\alpha^{-1} \Fnorm{G}_{\tilde{a}+2\alpha,\tilde{b}+\beta,m-1} \le
b \alpha^{-1} \Fnorm{G}_{\tilde{a}+2\alpha,\tilde{b}+\beta,m} \\
& \le &
b \alpha^{-1} \Fnorm{G}_{a,b,m} \le  \Fnorm{G}_{a,b,m} \le \bar{d}.
\end{eqnarray*}
We have used that
$\Phi(D_{\tilde{a}+\alpha,\tilde{b}})\subset
 D_{\tilde{a}+2\alpha,\tilde{b}+\beta}$ to bound
 $\Snorm{G\circ\Phi}_{\tilde{a},\tilde{b},m-1}$.

\item
If $\tilde{a} \le a-3\alpha$ and $\tilde{b} \le b - 2\beta$, then
\begin{eqnarray*}
\Fnorm{\tilde{G}_3}_{\tilde{a},\tilde{b},m} & \le &
M_3 \alpha^{-1}
\Fnorm{\Psi}_{\tilde{a}+3\alpha,\tilde{b}+2\beta,m-1}
\Fnorm{F \circ \Phi - A}_{\tilde{a}+\alpha,\tilde{b},m-1} \\
& \le &
M_3 \alpha^{-1} \Omega
\Fnorm{G^\bullet}_{\tilde{a}+3\alpha,\tilde{b}+2\beta,m}
\sigma \bar{d} \\
& \le &
\Omega \sigma M_3 \alpha^{-1} \bar{d} d^\bullet.
\end{eqnarray*}
The first inequality follows from Lemma~\ref{lem:Difference} with
$\Delta=\Psi$, $L=A$, $\Gamma_1=\mathbf{0}$,
$\Gamma_2=F\circ\Phi-A$, $l=m-1$, and $n=m-1$, so that
$M_3=M(\alpha,b,\beta,0,\sigma\bar{d}, m-1,m-1)$.

\item
If $\tilde{a} \le a-6\alpha$ and $\tilde{b} \le b - 4\beta$, then
\begin{eqnarray*}
\Fnorm{\tilde{G}_4}_{\tilde{a},\tilde{b},m} & \le &
\alpha^{-1} \Snorm{(\Upsilon+\Psi)\circ(F\circ\Phi)}_{\tilde{a}+\alpha,\tilde{b},m} \le
\alpha^{-1} \Snorm{\Upsilon + \Psi}_{\tilde{a}+3\alpha,\tilde{b}+2\beta,m} \\
& \le &
\alpha^{-1} \Fnorm{\Upsilon + \Psi}_{\tilde{a}+\alpha,\tilde{b},m} \le
M_4 \alpha^{-2} \Fnorm{\Psi}^2_{\tilde{a}+6\alpha,\tilde{b}+4\beta,m-1} \\
& \le &
M_4 \alpha^{-2} \left(\Omega\Fnorm{G^\bullet}_{\tilde{a}+6\alpha,\tilde{b}+4\beta,m}\right)^2 \\
& \le &
\Omega^2 M_4 \alpha^{-2} (d^\bullet)^2 \le
\Omega^2 M_4 \alpha^{-2} \bar{d} d^\bullet.
\end{eqnarray*}
The second inequality uses the inclusion
$(F\circ\Phi)(D_{\tilde{a}+\alpha,\tilde{b}})\subset
 D_{\tilde{a}+3\alpha,\tilde{b}+2\beta}$.
The fourth one follows from Lemma~\ref{lem:Inverse} with
$M_4=M(\alpha,b,\beta,\Omega\bar{d},\Omega\bar{d},m-1,m-1)$.
We need to verify the hypothesis $M_*bp<1$ in this last lemma.
It turns out that $M_*bp=M_4b\Omega\bar{d}=\Or(b^{m-2})$,
so it suffices to take $0<b\le\bar{b}$, with $\bar{b}$ small enough.

\item
If $\tilde{a} \le a-6\alpha$ and $\tilde{b} \le b - 4\beta$, then
\begin{eqnarray*}
\Fnorm{\tilde{G}^*}_{\tilde{a},\tilde{b},m} & \le &
\Fnorm{G^*}_{\tilde{a},\tilde{b},m} +
\textstyle \sum_{j=2}^4 \Fnorm{\tilde{G}_j}_{\tilde{a},\tilde{b},m} \le
d^* + \tilde{M} \bar{d} d^\bullet,\\
\Fnorm{\tilde{G}^\bullet}_{\tilde{a},\tilde{b},m} & \le &
\textstyle \sum_{j=1}^4 \Fnorm{\tilde{G}_j}_{\tilde{a},\tilde{b},m} \le
\left( \rme^{-2\pi K (a-\tilde{a})}  + \tilde{M} \bar{d} \right) d^\bullet,
\end{eqnarray*}
where $\tilde{M} = \Omega \alpha^{-1}\left(M_2+\sigma M_3+\Omega \alpha^{-1}M_4\right)$
and the constants $M_j$, $j=2,3,4$, have been defined previously.
\end{itemize}

\noindent
\emph{Step 4: Choice of the loss of analyticity domain.}
We set $\alpha = b$ and $\beta = \mu b^2/4$.
If $\bar{b}>0$ is small enough,
then conditions~(\ref{eq:Conditions}) hold for all $0<b\le\bar{b}$.
In addition,
\[
M=M(\alpha,b,\beta,c_1,c_2,l,n)=\Or(b^{l-2})\qquad
\mbox{as $b \to 0^+$},
\]
where $M$ is the expression introduced in Lemma~\ref{lem:Difference}.
If we take $\tilde{a}=a-6\alpha$ and $\tilde{b}=b-4\beta$,
then the bounds of the previous step imply that
\[
\Fnorm{\tilde{G}^*}_{\tilde{a},\tilde{b},m} \le
d^* + \tilde{M}\bar{d} d^\bullet,\qquad
\Fnorm{\tilde{G}^\bullet}_{\tilde{a},\tilde{b},m} \le
\left( \rme^{-12\pi K b} + \tilde{M}\bar{d} \right) d^\bullet,
\]
where $\tilde{M}=\tilde{M}(b;d,m,s)=
\Omega b^{-1}\left(M_2+\sigma M_3+\Omega b^{-1}M_4\right)=\Or(b^{m-5})$.
We recall that $m\ge 6$, $0 <\rho < s < 1$, and $K=s/b$.
Hence, if $0 < b \le \bar{b}$ and $\bar{b}$ is small enough, then
\[
\Fnorm{\tilde{G}^*}_{\tilde{a},\tilde{b},m} \le
d^* + \rme^{-12\pi \rho} d^\bullet,\qquad
\Fnorm{\tilde{G}^\bullet}_{\tilde{a},\tilde{b},m} \le
\rme^{-12 \pi \rho} d^\bullet.
\]
Indeed, $\tilde{M} \bar{d} \le
\rme^{-12 \pi \rho} - \rme^{-12 \pi s} \le
\rme^{-12 \pi \rho}$ if we take a small enough value of $\bar{b}$.

This ends the proof of the iterative lemma.

\subsection{Proof of Theorem~\ref{thm:Technical}}\label{sec:ProofThmTechnical}

Set $\rho = \sqrt{r} \in (0,1)$, $\mu = 6(1-\rho)/(\bar{a}-\breve{a})$,
and $\Omega = \Omega(\sqrt{\rho})$,
where the function $\Omega(s)$ is defined in Lemma~\ref{lem:LinearDouble}.
Let $\bar{b}$ be the positive constant associated
to the integer $m \ge 6$ in Lemma~\ref{lem:Iterative},
the numbers $\bar{a},\bar{d},\mu>0$, and the exponent $\rho\in(0,1)$.
Let  $c_1 = c_1(r) = \Omega \sum_{n \ge 0} \rme^{-12\pi \rho n}$,
$c_2 = c_2(r) = \sum_{n \ge 1} \rme^{-12\pi \rho n}$, and
$c_3 = c_3(r) = \rme^{12\pi \rho}$.

Let us check that $\bar{b}$, $c_1$, $c_2$, and $c_3$ satisfy
the properties given in Theorem~\ref{thm:Technical}.

Let $a_0 = \bar{a}$, $d_0^* = \bar{d}^*$, $d_0^\bullet = \bar{d}^\bullet$,
$0 < b_0 = \breve{b}/\rho \le \bar{b}$,
$F_0 = \bar{F} =  A + \bar{G}$ be the map given in~(\ref{eq:barFv}),
$G_0^* = \pi^*(\bar{G})$, and $G_0^\bullet = \pi^\bullet(\bar{G})$.
By recursively applying Lemma~\ref{lem:Iterative},
we obtain a sequence of changes of variables $\Phi_n = \Identity + \Psi_n$,
with $\Psi_n \in \mathcal{X}_{a_{n-1},b_{n-1},m-1}$,
and a sequence of maps $F_n = A + G_n$,
with $G_n = G_n^* + G_n^\bullet$,
$G_n^* \in \mathcal{X}_{a_n,b_n,m}^*$ and $G_n^\bullet \in \mathcal{X}_{a_n,b_n,m}^\bullet$,
such that
\[
\Fnorm{\Psi_n}_{a_{n-1},b_{n-1},m-1} \le \Omega d_{n-1}^\bullet,\quad
\Fnorm{G_n^*}_{a_n,b_n,m}\le d_n^*,\quad
\Fnorm{G_n^\bullet}_{a_n,b_n,m}\le d_n^\bullet,
\]
with $a_{n+1} = a_n-6b_n$,
$b_{n+1} = b_n - \mu b_n^2$,
$d_{n+1}^* = d_n^* + \rme^{-12\pi \rho} d_n^\bullet$, and
$d_{n+1}^\bullet = \rme^{-12\pi \rho} d_n^\bullet$.

Let $N$ be the biggest integer satisfying $N b_0 \le (\bar{a}-\breve{a})/6$.
The sequences $(a_n)_{0 \le n \le N}$, $(b_n)_{0 \le n \le N}$,
and $(d_n^\bullet)_{0 \le n \le N}$ are decreasing.
The sequence $(d_n^*)_{0 \le n \le N}$ is increasing.
Indeed,
\begin{eqnarray*}
a_N & = & a_{N-1}-6b_{N-1} \ge a_{N-1}-6b_0 \ge \cdots \ge a_0 - 6Nb_0 \ge \breve{a}, \\
b_N & = & b_{N-1} - \mu b^2_{N-1} \ge b_{N-1} - \mu b^2_0 \ge
\cdots \ge (1- \mu N b_0)b_0 \ge \rho b_0 = \breve{b},\\
d_N^\bullet & \le & \rme^{-12 \pi \rho} d_{N-1}^\bullet \le
\cdots \le \rme^{-12 \pi \rho N} d_0^\bullet \le
c_3 \rme^{-2\pi r (\bar{a}-\breve{a})/\breve{b}} \bar{d}^\bullet,\\
d_N^* & \le & d_{N-1}^* + \rme^{-12 \pi \rho} d_{N-1}^\bullet \le \cdots \le
d_0^* + \left( \textstyle \sum_{n=1}^N \rme^{-12\pi \rho n} \right)
d_0^\bullet \le
\bar{d}^* + c_2 \bar{d}^\bullet,
\end{eqnarray*}
and $d_n^* + d_n^\bullet \le d_N^* + d_0^\bullet \le
\bar{d}^* + (1+c_2) \bar{d}^\bullet \le \bar{d}$ for all $n = 0,\ldots,N$.

We can apply $N$ times the iterative lemma.
Let $\breve{F} = A + \breve{G} = A + G_N = F_N$ be the map obtained
after those $N$ steps. Then
\begin{eqnarray*}
\Fnorm{\pi^*(\breve{G})}_{\breve{a},\breve{b},m} & \le &
\Fnorm{\pi^*(G_N)}_{a_N,b_N,m} \le
d_N^* \le \bar{d}^* + c_2 \bar{d}^\bullet, \\
\Fnorm{\pi^\bullet(\breve{G})}_{\breve{a},\breve{b},m} & \le &
\Fnorm{\pi^\bullet(G_N)}_{a_N,b_N,m} \le
d_N^\bullet \le c_3 \rme^{-2 \pi r (\bar{a}-\breve{a})/\breve{b}} \bar{d}^\bullet.
\end{eqnarray*}
Finally, let $\breve{\Phi} = \Phi_N \circ \cdots \circ \Phi_1$ be the change
of variables such that $\breve{F} = \breve{\Phi}^{-1} \circ \bar{F} \circ \breve{\Phi}$.
We want to check that $\breve{\Phi} = \Identity + \breve{\Psi}$ for some
$\breve{\Psi} \in \mathcal{X}_{\breve{a},\breve{b},m-1}$
such that $\Snorm{\breve{\Psi}}_{\breve{a},\breve{b},m-1} \le c_1 \bar{d}^\bullet$.
We note that
\[
\breve{\Psi} = \Psi_1 + \cdots + \Psi_N,
\]
where each term of the above summation is evaluated at a different argument.
Nevertheless, those arguments are not important when computing the sup-norm:
\[
\Snorm{\breve{\Psi}}_{\breve{a},\breve{b},m-1} \le
\sum_{n=1}^N \Snorm{\Psi_n}_{a_{n-1},b_{n-1},m-1} \le
\Omega \sum_{n=0}^{N-1} d_n^\bullet \le
\Omega d_0^\bullet \sum_{n=0}^{N-1} \rme^{-12 \pi \rho n} =
c_1 \bar{d}^\bullet.
\]

This ends the proof of Theorem~\ref{thm:Technical}.

\subsection{Proof of Theorem~\ref{thm:Remainder}}\label{sec:ProofThmRemainder}

Let us begin with a simple, but essential, chain of inequalities associated
to certain analyticity strip widths that will appear along the proof.
If $\alpha \in (0,a_*)$, then there exists $r \in (0,1)$, $\bar{b} > 0$,
and some analyticity strip widths $a_2$, $\bar{a} = a_m$, and $\breve{a} = \bar{b}$,
such that
\begin{equation}\label{eq:achain}
0 < \bar{b} =: \breve{a} := \bar{a}-(1+\bar{b})\alpha/r <
\bar{a} := a_m < a_2 < a_*.
\end{equation}
The first two reductions (that is, from $a_*$ to $a_2$ and from $a_2$ to $a_m$)
are as small as we want.
The third reduction (from $\bar{a} = a_m$ to $\breve{a} = \bar{a}-(1+\bar{b})\alpha/r$)
should be a little bigger than $\alpha$ in order to get the desired exponentially
small upper bound with the exponent $\alpha$.
The fourth reduction (that is, from $\breve{a} = \bar{b}$ to $0$) is also small,
since $\bar{b}$ can be taken as small as necessary.

This decreasing positive sequence of analyticity strip widths
is associated to a similar sequence of analyticity radii.
To be precise, we will construct a sequence of the form
\[
b < \breve{b} < \bar{b} \le b_m < b_2 < b_*,\qquad
\breve{b} := b + b^2 < \bar{b} \sqrt{r}.
\]
The inequality~$\bar{b} \le b_m$ does not correspond to a true reduction,
but to a restriction on the size of $\bar{b}$.
Note that we have consumed all the analyticity strip width
after the last reduction, but we still keep a positive analyticity
radius $b$.

We split the proof in the eight steps.

\noindent
\emph{Step 1: Control of the Fourier norm.}
If the analytic map $f$ satisfies the properties (i)--(iii)
listed in Lemma~\ref{lem:Order2},
then the map $F_2 = A + G_2 := f$ is real analytic and
has the intersection property on the cylinder $\Tset \times (-b_*,b_*)$,
can be extended to the complex domain $D_{a_*,b_*}$,
and has the form~(\ref{eq:Order2}).
The Fourier norm $\Fnorm{G_2}_{a_*,b_*,2}$ may be infinite,
but $\Fnorm{G_2}_{a_2,b_2,2} < \infty$ for any $a_2 \in (0,a_*)$
and $b_2 \in (0,b_*)$.

\noindent
\emph{Step 2: Application of the averaging lemma.}
Once fixed an integer $m \ge 6$ and any $a_m \in (0,a_2)$,
we know from Lemma~\ref{lem:BeyondAllOrders}
that there exist an analytical radius $b_m \in (0,b_2)$ and
a change of variables of the form $\Phi_m = \Identity + \Psi_m$
for some $\Psi_m \in \mathcal{X}_{a_m,b_m,1}$
such that the transformed map
$F_m = \Phi_m^{-1} \circ F_2 \circ \Phi_m$
is real analytic, has the intersection property on the cylinder
$\Tset \times (-b_m,b_m)$, and has the form
$F_m = A + G_m$ for some $G_m \in \mathcal{X}_{a_m,b_m,m}$.

\noindent
\emph{Step 3: Application of Theorem~\ref{thm:Technical}.}
Let $r \in (0,1)$ be the number that appears in~(\ref{eq:achain}).
Set $\bar{F} = A + \bar{G} = A + G_m$, $\bar{a} = a_m$, and
$\bar{d} =
 \Fnorm{\pi^*(\bar{G})}_{\bar{a},b_m,m} +
 (1+c_2(r)) \Fnorm{\pi^\bullet(\bar{G})}_{\bar{a},b_m,m}$.

Let $\bar{b} = \bar{b}(m,\bar{a},\bar{d},r) > 0$ be the constant
stated in Theorem~\ref{thm:Technical}.
We can assume that $\bar{b} \le b_m$ and the condition~(\ref{eq:achain}) holds,
by taking a smaller $\bar{b} > 0$ if necessary.
Let $b'_* \in (0,b_*)$ be defined by $b'_* + (b'_*)^2 = \bar{b} \sqrt{r}$.
Fix any $b \in (0,b'_*)$.
Set $\breve{a} = \bar{a}-(1+\bar{b})\alpha/r$ and
$\breve{b} = b + b^2 \le b'_* + (b'_*)^2 = \bar{b} \sqrt{r}$.

If $\bar{d}^*$ and $\bar{d}^\bullet$ are the norms defined in~(\ref{eq:barFv}),
then $\bar{d}^* + (1+c_2 \bar{d}^\bullet) \le \bar{d}$.
Hence, we can apply Theorem~\ref{thm:Technical} to obtain a change of variables
$\breve{\Phi} = \Identity + \breve{\Psi}$,
with $\breve{\Psi} \in \mathcal{X}_{\breve{a},\breve{b},m}$ and
$\Snorm{\breve{\Psi}}_{\breve{a},\breve{b},m-1} \le
c_1 \bar{d}^\bullet \le c_1 \bar{d}$,
and a transformed map
$\breve{F} =  A + \breve{G} = \breve{\Phi}^{-1} \circ \bar{F} \circ \breve{\Phi}$,
with $\breve{G} \in \mathcal{X}_{\breve{a},\breve{b},m}$,
$\Snorm{\pi^\bullet(\breve{G})}_{\breve{a},\breve{b},m} \le
\Fnorm{\pi^\bullet(\breve{G})}_{\breve{a},\breve{b},m} \le
c_3 \rme^{-2\pi r (\bar{a}-\breve{a})/\breve{b}} \bar{d}^{\bullet} \le
c_3 \rme^{-2\pi \alpha (1+\bar{b})/\breve{b}} \bar{d} \le
c_3 \rme^{-2 \pi \alpha/b} \bar{d}$, and
$\Snorm{\pi^*(\breve{G})}_{\breve{a},\breve{b},m} \le
 \Fnorm{\pi^*(\breve{G})}_{\breve{a},\breve{b},m} \le
 \bar{d}$.

\noindent
\emph{Step 4: Uniform estimates on the change $\Phi = \Phi_m \circ \breve{\Phi}$}.
By construction, $\breve{\Phi} = \Identity + \breve{\Psi}$,
with $\breve{\Psi} \in \mathcal{X}_{\breve{a},\breve{b},m-1}$ and
$\Snorm{\breve{\Psi}}_{\breve{a},\breve{b},m-1} \le \breve{M}$,
where the constant $\breve{M} := c_1 \bar{d}$ does not depend on $b$.
Thus,
\[
\breve{\Phi}(x,y) =
\left(x + y^{m-1} \breve{\psi}_1(x,y), y + y^m \breve{\psi}_2(x,y) \right)
\]
for some functions $\breve{\psi}_j(x,y)$ analytic
on $D_{\breve{a},\breve{b}} = D_{\bar{b},b+b^2}$ such that
$|\breve{\psi}_j|_{\bar{b},b+b^2} \le \breve{M}$.
The Cauchy estimates imply that
\[
|\breve{\psi}_j(x,y)| \le \breve{M},\qquad
|\partial_1 \breve{\psi}_j(x,y)| \le \bar{b}^{-1} \breve{M},\qquad
|\partial_2 \breve{\psi}_j(x,y)| \le b^{-2} \breve{M},
\]
for all $(x,y) \in \Tset \times B_b$ and, in particular,
for all $(x,y) \in \Tset \times (-b,b)$.
Hence,
\[
\breve{\Phi}(x,y) = \left( x+\Or(y^{m-1}),y+\Or(y^m) \right),\qquad
\det\left[ \breve{\Phi}(x,y) \right] = 1 + \Or(y^{m-2}),
\]
for all $(x,y) \in \Tset \times (-b,b)$, where the
$\Or(y^{m-2})$, $\Or(y^{m-1})$, and $\Or(y^m)$ terms are  uniform in $b$.
We recall that $m \ge 6$ and the change $\Phi_m$ satisfies
properties~(\ref{eq:ChangePhi_m}),
so the complete change $\Phi = \Phi_m \circ \breve{\Phi}$
satisfies the properties stated in Theorem~\ref{thm:Remainder}.

\noindent
\emph{Step 5: Exponentially small bound on the remainder $G$.}
After all these changes of variables,
we have the map $F = A + G := \breve{F}$,
with $G := \breve{G} \in
\mathcal{X}_{\breve{a},\breve{b},m}$,
$\Snorm{\pi^*(G)}_{\breve{a},\breve{b},m} \le \bar{d}$
and $\Snorm{\pi^\bullet(G)}_{\breve{a},\breve{b},m} \le
c_3 \bar{d} \rme^{-2 \pi \alpha/b}$.
We can bound $G^* = \pi^*(G)$ by using the bound on
$G^\bullet = \pi^\bullet(G)$ and the intersection property
of $F$ on the cylinder $\Tset \times (-b,b)$.
We recall that if
\[
G(\xi,\eta) = (\eta^m g_1(\xi,\eta), \eta^{m+1} g_2(\xi,\eta) ),
\]
for some $g_1,g_2 \in \mathcal{X}_{a,b}$, then
\[
G^*(\xi,\eta) = (0, \eta^{m+1} g_2^*(\eta) ), \qquad
G^\bullet(\xi,\eta) =
(\eta^m g_1(\xi,\eta), \eta^{m+1} g^\bullet_2(\xi,\eta) ),
\]
where $g_2^*(\eta)$ is the average of $g_2(\xi,\eta)$
and $g_2^\bullet = g_2 - g_2^*$.
Fix any $\eta_0 \in (-b,b)$.
We know that
\[
F(\Tset \times \{ \eta_0 \}) \cap \big( \Tset \times \{ \eta_0 \} \big)
\neq \emptyset.
\]
Therefore, there exists~$\xi_0 \in \Tset$ such that
$g_2^*(\eta_0) + g_2^\bullet(\xi_0,\eta_0) = 0$,
and so
\[
|g_2^*(\eta)| \le
\sup_{\Tset \times B_b} |g_2^\bullet| \le
\Snorm{G^\bullet}_{\breve{a},\breve{b},m} \le
c_3 \bar{d} \rme^{-2 \pi \alpha/b},\qquad
\forall \eta \in (-b,b).
\]
This implies that
$|g_j(\xi,\eta)| \le \Snorm{G}_{\breve{a},\breve{b},m} \le
 \Snorm{G^*}_{\breve{a},\breve{b},m} +
 \Snorm{G^\bullet}_{\breve{a},\breve{b},m} \le
 2c_3 \bar{d} \rme^{-2 \pi \alpha/b}$
for all $(\xi,\eta) \in \Tset \times (-b,b)$.

\noindent
\emph{Step 6: Exponentially small bounds on some derivatives of the remainder.}
We recall that
\[
\max \left \{ |g_1(\xi,\eta)|,|g_2^\bullet(\xi,\eta)| \right\} \le
\Snorm{G^\bullet}_{\breve{a},\breve{b},m} \le
c_3 \bar{d} \rme^{-2 \pi \alpha/b}
\]
for all $(\xi,\eta) \in D_{\breve{a},\breve{b}} = D_{\bar{b},b+b^2}$.
Thus, we get from $\partial_1 g_2 = \partial_1 g_2^\bullet$ and
the Cauchy estimates that
\begin{eqnarray*}
|\partial_1 g_j(\xi,\eta)| & \le &
c_3 \bar{d} \ \bar{b}^{-1} \rme^{-2 \pi \alpha/b},\\
\max \left \{ |\partial_2 g_1(\xi,\eta)|,|\partial_2 g_2^\bullet(\xi,\eta)| \right\}
& \le & c_3 \bar{d} b^{-2} \rme^{-2 \pi \alpha/b},
\end{eqnarray*}
for all $(\xi,\eta) \in \Tset \times (-b,b)$.

\noindent
\emph{Step 7: A crude bound on the derivative of $g_2^*$.}
We recall that $G^*(\xi,\eta) = (0,\eta^{m+1}g_2^*(\eta))$, so
\[
|g^*_2(\eta)| \le
\Snorm{G^*}_{\breve{a},\breve{b},m} \le \bar{d},\qquad
\forall \eta \in B_{\breve{b}} = B_{b+b^2}.
\]
Therefore, the Cauchy estimates imply that
$|(g_2^*)'(\eta)| \le b^{-2} \bar{d}$ for all $\eta \in B_b$ and,
in particular, for all $\eta \in (-b,b)$.

\noindent
\emph{Step 8: Computation of the constant $K$.}
By combining the inequalities obtained in Steps~5--7, we get that
$|g_j(\xi,\eta)| \le K \rme^{-2 \pi \alpha/b}$ and
$|\partial_i g_j(\xi,\eta)| \le K b^{-2}$ for all $\Tset \times (-b,b)$,
provided
\[
K = \bar{d} \max \left\{ 2c_3, c_3 + 1 \right\}.
\]

This ends the proof of Theorem~\ref{thm:Remainder}.

\section{Proof of Theorem~\ref{thm:Asymptotic}}\label{sec:ProofThmAsymptotic}

\subsection{A space of matrix functions}

Henceforth, let $I_b = (-b,b) \subset \Rset$ and
$S_b = \Tset \times I_b$ with $b>0$.
Let $S$ be any compact subset of $S_b$.
Let $\mu \in \Nset$.
Let $\mathcal{M}_{S,\mu}$ be the set of all matrix functions
$\Gamma : S \to \mathcal{M}_{2\times 2}(\Rset)$ of the form
\[
\Gamma(\xi,\eta) =
\left(
\begin{array}{cc}
\eta^\mu \gamma_{11}(\xi,\eta) & \eta^{\mu-1} \gamma_{12}(\xi,\eta) \\
\eta^{\mu+1} \gamma_{12}(\xi,\eta) & \eta^\mu \gamma_{22}(\xi,\eta)
\end{array}
\right),
\]
for some continuous functions $\gamma_{ij}:S \to \Rset$.
The set $\mathcal{M}_{S,\mu}$ is a Banach space with the norm
\[
\| \Gamma \|_{S,\mu} =
\max \left\{ |\gamma_{ij}(\xi,\eta)| :
             (\xi,\eta) \in S, \ 1\le i,j \le 2 \right\}.
\]

\begin{lem}\label{lem:BoundGammas}
Let $S \subset S_b$, $\Gamma \in \mathcal{M}_{S,\mu}$,
$\Delta \in \mathcal{M}_{S,\nu}$, and
$A = \left( \begin{array}{cc} 1 & 1 \\ 0 & 1 \end{array}\right)$.
Let $k \in \Nset$ and $j \in \Zset$.
Let $f:S \to S$ be a map of the form~(\ref{eq:NeishtadtForm})
with $|g_2(\xi,\eta)|\le K_0$ for all $(\xi,\eta) \in S$.
Then:
\begin{enumerate}
\item
$\Gamma\Delta \in \mathcal{M}_{S,\mu+\nu}$ and
$\| \Gamma\Delta \|_{S,\mu+\nu} \le
 2 \| \Gamma \|_{S,\mu} \| \Delta \|_{S,\nu}$.
\item
$A^k \Gamma\in \mathcal{M}_{S,\mu}$ and
$\| A^k \Gamma \|_{S,\mu} \le (1+bk) \| \Gamma \|_{S,\mu}$.
\item
$\Gamma A^k\in \mathcal{M}_{S,\mu}$ and
$\| \Gamma A^k \|_{S,\mu} \le (1+bk) \| \Gamma \|_{S,\mu}$.
\item
$\Gamma \circ f^j \in \mathcal{M}_{S,\mu}$ and
$\| \Gamma \circ f^j \|_{S,\mu} \le (1+K_0 b^m)^{(\mu+1)|j|} \| \Gamma \|_{S,\mu}$.
\end{enumerate}
\end{lem}

\proof
It is a straightforward computation.
\qed

\subsection{A technical lemma}

Let $f$ be an analytic map of the form~(\ref{eq:NeishtadtForm}).
Let $p$ and $q$ be two relatively prime integers.
There exist two curves $R = \Graph \zeta$ and $\hat{R} = \Graph \hat{\zeta}$
and two RICs $R_\pm = \Graph \zeta_\pm$ with Diophantine rotation numbers
$\omega_- < p/q < \omega_+$,
all four contained in a small neighborhood of $\Tset \times \{ p/q \}$,
such that $f^q$ projects $R$ onto $\hat{R}$ along the vertical direction
and $R$ and $\hat{R}$ are contained in the strip of the cylinder
enclosed by the RICs $R_\pm$.
Following Birkhoff~\cite[Section VI]{Birkhoff1966} and
Arnold~\cite[Section 20]{ArnoldAvez1968},
all $(p,q)$-periodic points of $f$ are contained in $R \cap \hat{R}$.
Besides, we will see later on that the (geometric) area enclosed between $R$
and $\hat{R}$ is an upper bound of the quantities $\Delta^{(p,q)}$.
These are the reasons for the study of $R$ and $\hat{R}$.

Let us prove that these four curves exist for big enough periods $q$.
In this case, ``big enough'' only depends on the size of the
nonintegrable terms of~$f$,
the size of the neighborhood of $\Tset \times \{ p/q \}$, the exponent $m$,
and the winding number $p$.
On the contrary, it \emph{does not depend} on the particular map at hand.
Therefore, every time that we ask $q$ to be ``big enough'' along the proof
of the following lemma, it only depends on the quantities $K_0>0$, $c>1$,
$m \ge 4$, $p \in \Zset \setminus \{ 0 \}$, and $q'_* \in \Nset$ fixed at the first
line of the next statement.

\begin{lem}\label{lem:Preliminary}
Let $K_0 > 0$,  $c>1$, $m \ge 4$, and $p,q'_* \in \Nset$.
Let $q \ge q'_*$ be an integer relatively prime with $p$.
Set $b = c^2 p/q$.
Let $f: S_b \to \Tset \times \Rset$ be an analytic map of the
form~(\ref{eq:NeishtadtForm}) such that $|g_j(\xi,\eta)| \le K_0$
and $|\partial_i g_j(\xi,\eta)| \le K_0 b^{-2}$ for all $(\xi,\eta) \in S_b$.
Let $(\xi_q,\eta_q) = f^q(\xi,\eta)$.
Let $I = (p/c^2 q,c^2 p/q)$, $I_- = (p/c^2 q,p/c q)$, and $I_+ = (cp/q,c^2 p/q)$.
There exists $q''_* = q''_*(K_0,c, m, p, q'_*) \ge q'_*$ such that, if $q \ge q''_*$,
the following properties hold:
\begin{enumerate}
\item
The map $f$ has two RICs $R_\pm \subset \Tset \times I_\pm \subset S_b$
whose internal dynamics is conjugated to a rigid rotation
of angles $\omega_\pm \in I_\pm$, respectively;
\item
If $S$ is the compact subset of $S_b$ enclosed by $R_-$ and $R_+$, then
\begin{equation}\label{eq:TwistCondition}
\frac{\partial \xi_q}{\partial \eta}(\xi,\eta) > 0,\qquad
\forall (\xi,\eta) \in S;
\end{equation}
\item \label{item:Zeta}
There exist two unique analytic functions
$\zeta: \Tset \to I$ and $\hat{\zeta}: \Tset \to I$ such that
\begin{equation}\label{eq:VerticalCurves}
f^q (\xi,\zeta(\xi)) = (\xi,\hat{\zeta}(\xi)), \qquad \forall \xi \in \Tset,
\end{equation}
and all the $(p,q)$-periodic points of the restriction $f|_S$
are contained in $\Graph \zeta$.
\end{enumerate}
The same statement holds if $p$ is a negative integer,
$b = c^2|p|/q$, $I = (c^2 p/q, p/c^2 q)$, $I_- = (c^2 p/q,cp/q)$,
and $I_+ = (p/cq,p/c^2 q)$.
\end{lem}

\proof
Let us assume $p>0$. The case $p<0$ is analogous.

First, the existence of the RICs $R_-$ and $R_+$ follows from some quantitative
estimates in KAM theory established by Lazutkin~\cite[Theorem 2]{Lazutkin1973}.
To be precise, Lazutkin proved that there exists $b'_* = b'_*(K_0) > 0$
such that if $\omega \in (-b'_*,b'_*)$ satisfies the Diophantine condition
\begin{equation}\label{eq:Diophantine}
2 |\omega - i/j| \ge |i| j^{-4}
\end{equation}
for all integers $j \ge 1$ and $i$, then $f$ has a RIC
$R = \{ \eta = \omega + \Or(\omega^m) \}$ whose internal dynamics
is $C^l$-conjugated to a rigid rotation of angle $\omega$,
for a suitable $l \ge 1$. The conjugation is $\Or(1/q^{m-1})$-close to the identity.
Item~(i) follows directly from this estimate,
because there exist some real numbers
 $\omega_+\in (c^{4/3}p/q,c^{5/3}p/q )\subset I_+$
and $\omega_-\in (p/c^{5/3} q,p/c^{4/3}q )\subset I_-$
 satisfying the Diophantine
condition~(\ref{eq:Diophantine}),
provided $q$ is big enough.

Second, let us check that the power map $f^q$ satisfies~(\ref{eq:TwistCondition}).
The compact subset $S \subset S_b$ is invariant by $f$,
because it is delimited by RICs.
Thus, all powers $(\xi_j,\eta_j) = f^j(\xi,\eta)$ are well defined on $S$.
We write $\rmD f (\xi,\eta) = A + \Gamma(\xi,\eta)$,
where $A$ was introduced in Lemma~\ref{lem:BoundGammas}.
Next, we compute the differential of the power map:
\begin{equation}\label{eq:Dfq}
\rmD f^q =
(A+\Gamma_q)\cdots (A+\Gamma_1) =
A^q + \Delta_1 + \cdots + \Delta_q,
\end{equation}
where $\Gamma_j = \Gamma \circ f^j$ and
$\Delta_l$ is the sum of all the products of the form
$A^{k_1} \Gamma_{j_1} \cdots A^{k_l} \Gamma_{j_l} A^{k_{l+1}}$,
with $k_i \ge 0$, $q \ge j_1 > j_2 > \cdots > j_l \ge 1$,
and $q= l+\sum_{i=1}^{l+1} k_i$.
These products are elements of $\mathcal{M}_{S,m l}$,
because $\Gamma \in \mathcal{M}_{S,m}$.
Indeed, if $C = \| \Gamma \|_{S,m}$, then
\begin{eqnarray*}\fl
\left\| A^{k_1} \Gamma_{j_1} \cdots
         A^{k_l} \Gamma_{j_l} A^{k_{l+1}} \right\|_{S,m l}
& \le &
2^{l-1} \|A^{k_1} \Gamma_{j_1}\|_{S,m}\cdots
       \|A^{k_{l-1}} \Gamma_{j_{l-1}}\|_{S,m}
       \|A^{k_l} \Gamma_{j_l}A^{k_{l+1}}\|_{S,m} \\
& \le &
2^{l-1} C^l \textstyle
\prod_{i=1}^{l+1} (1 + b k_i)
\prod_{i=1}^l (1+K_0 b^m)^{(m+1) j_i} \\
& \le &
2^{l-1} C^l \exp \left(
\textstyle \sum_{i=1}^{l+1} b k_i +
                                   \sum_{i=1}^l (m+1) K_0 b^m j_i \right) \\
&\le &
2^{l-1} C^l \rme^{c^2 p + (m+1)K_0c^{2 m} p^m} =
C'(2C)^l,
\end{eqnarray*}
where we have used Lemma~\ref{lem:BoundGammas},
inequality $1+x\le \rme^x$ for $x\ge 0$,
$\sum_{i=1}^{l+1} k_i\le q$, $l \le q$,
$j_i \le q$, $b = c^2 p/q$, and $m \ge 4$.
We have also defined $C' = \rme^{c^2 p + (m+1) K_0 c^{2 m} p^m} /2$.

The matrix $\Delta_l$ is the sum of the products
with precisely $l$ factors $\Gamma_j$.
This shows that there are ${q \choose l}$ terms inside $\Delta_l$.
Therefore, $\Delta_l\in\mathcal{M}_{S,m l}$ and
\begin{equation}\label{eq:BoundDelta}
\|\Delta_l\|_{S,m l} \le
{q \choose l}
\left\| A^{k_1} \Gamma_{j_1} \cdots
         A^{k_l} \Gamma_{j_l} A^{k_{l+1}} \right\|_{S,m l} \le
    C'(2Cq)^l.
\end{equation}
The element of the first row and second column of $A^q$
is equal to $q$, so
\begin{eqnarray*}
\left|\frac{\partial \xi_q }{\partial \eta}(\xi,\eta) - q \right|
& \le &
C'\sum_{l=1}^q (2 C q)^l b^{m l-1} \le
\frac{C'}{b} \sum_{l=1}^q \left(2 K_0 q b^{m-2} \right)^l \\
& \le &
4 C' K_0 q b^{m-3} \le
4 C' K_0 c^{2(m-3)} p^{m-3},
\end{eqnarray*}
for all $(\xi,\eta) \in S \subset S_b$,
which implies the twist condition~(\ref{eq:TwistCondition})
provided that  $q$ is big enough.
Here, we have used relation~(\ref{eq:Dfq}), bound~(\ref{eq:BoundDelta}),
$b = c^2 p/q$, and $m\ge 4$.
We have also used that $C \le K_0 b^{-2}$ and $2 K_0 q b^{m-2} \le 1/2$,
provided $q$ is big enough.

Third, we establish the existence of the functions
$\zeta,\hat{\zeta}:\Tset \to I$.
We know from Lazutkin~\cite{Lazutkin1973} that
$R_\pm = \Graph \zeta_\pm$ for some differentiable functions
$\zeta_\pm: \Tset \to I_\pm$.
We work with the lifts $F$, $\Xi_q$, and $Z_\pm$
of the objects $f$, $\xi_q$, and $\zeta_\pm$.
The RICs are invariant,
so $F^q(\xi,Z_\pm(\xi)) = (\Xi_\pm(\xi),Z_\pm(\Xi_\pm(\xi)))$ for some
differentiable functions $\Xi_\pm: \Rset \to \Rset$.
If we prove that there exist two unique analytic 1-periodic functions
$Z,\hat{Z}:\Rset\to I$ such that
\begin{equation}\label{eq:VerticalCurvesLift}
 F^q(\xi,Z(\xi)) = (\xi+p,\hat{Z}(\xi+p)),\qquad \forall \xi\in \Rset,
\end{equation}
then item~(iii) follows.
Since the dynamics of $F^q$ on $R_\pm$ is $C^l$-conjugated to
a rigid rotation of angle $q\omega_+$ through a $\Or(1/q^{m-1})$-close to the identity
conjugation,
\[
\Xi_+(\xi)= \xi+q\omega_+ +\Or(1/q^{m-1}) \ge \xi+cp+\Or(1/q^{m-1})>\xi+p
\]
provided that $q$ is big enough.
Analogously, we obtain $\Xi_-(\xi)<\xi+p$.
That is,
\[
\Xi_q(\xi,Z_-(\xi)) = \Xi_-(\xi) < \xi+p < \Xi_+(\xi) = \Xi_q(\xi,Z_+(\xi)),
\qquad \forall \xi\in \Rset.
\]
Since $\Xi_q(\xi,\eta)$ is analytic and strictly increasing for
$\eta \in (Z_-(\xi),Z_+(\xi)) \subset I$,
we deduce that there exists a unique function $Z: \Rset \to I$
such that $\Xi_q(\xi,Z(\xi)) = \xi + p$.

The function $S \ni (\xi,\eta) \mapsto G(\xi,\eta) := \Xi_q(\xi,\eta)-\xi-p$ is analytic
and $\frac{\partial G}{\partial \eta}(\xi,\eta) > 0$,
so $Z$ is analytic by the Implicit Function Theorem.
The $1$-periodicity of $Z$ follows from the uniqueness and the property
$F^q(\xi+1,\eta) = F^q(\xi,\eta)+(1,0)$.
Function $\hat{Z}: \Rset \to I$ is defined by means of
relation~(\ref{eq:VerticalCurvesLift}).
Finally, functions $\zeta, \hat{\zeta}:\Tset \to I$ are
the projections of $Z,\hat{Z}: \Rset \to I$.
\qed

\subsection{Proof of Theorem~\ref{thm:Asymptotic}: Case $(m,n)=(0,1)$}

If $(m,n) = (0,1)$, by hypothesis,
the map $g:\Tset\times I\to \Tset\times I$,
$(s,r)\mapsto(s_1,r_1)$, is an analytic exact twist map with
a $(a_*,b_*)$-analytic $(0,1)$-resonant RIC,
such that $\varrho_-\le 0 \le \varrho_+$.
The map $f = g^n = g$ satisfies the properties~(i)--(iii)
listed in Lemma~\ref{lem:Order2}
in some suitable coordinates $(x,y)$.
Let $(s,r) = \tilde{\Phi}(x,y)$ be the associated change of variables.
Let $\tilde{f} = \tilde{\Phi}^{-1} \circ f\circ \tilde{\Phi}$
be the new map  defined in the domain~(\ref{eq:Dab}).
Note that the $(a_*,b_*)$-analytic $(0,1)$-resonant RIC is
$C\equiv \{y=0\}$ in the $(x,y)$ coordinates.

Let $p$ be an integer such that $1 \le |p| \le L$.
Let $c \in (1,2)$ such that $\alpha<c^2\alpha <a_*$.
We take $c^2\alpha $ as the $\alpha$ appearing in Theorem~\ref{thm:Remainder},
$m=4$, and $b=c^2|p|/q$, provided that~$q$ is relatively prime with $p$
and is large enough so that $c^2|p|/q<b'_*=b'_*(\alpha)$.
That is, $q>q'_*:=c^2|p|/b'_*$.

Hence, there exist $K_0,K_1>0$, both independent of $q$,
and a change of coordinates $(x,y)=\Phi(\xi,\eta)$ such that
$\bar{f} = \Phi^{-1} \circ \tilde{f}\circ \Phi: S_b \to \Tset\times\Rset$
is an analytic map of the form~(\ref{eq:NeishtadtForm}) such that
$|g_j(\xi,\eta)| \le
 K_0 \rme^{-2\pi c^2 \alpha/b} =
 K_0 \rme^{-2\pi \alpha q/|p|} \le K_0$,
$|\partial_i g_j(\xi,\eta)| \le  K_0 b^{-2}$, and
$\sup \{ \left| \det[\rmD\Phi(\xi,\eta)] \right| \} \le  K_1$
for all $(\xi,\eta) \in S_b$.

The map $\bar{f}:S_b \to \Tset\times\Rset$
satisfies the hypotheses of Lemma~\ref{lem:Preliminary} for any $q\ge q'_*$.
Let $q_*$ be the maximum value of $q''_*$ among the integers $0 < |p|\le L$.
Let $R_\pm$ be the RICs with rotation numbers $\omega_\pm$ given
in Lemma~\ref{lem:Preliminary}.
Let $S$ be the compact subset of $S_b$ enclosed by $R_-$ and $R_+$.
Since $f$ is globally twist and $\varrho_- < \omega_- < p/q < \omega_+ < \varrho_+$,
all the Birkhoff $(p,q)$-periodic orbits of $f$ are
contained in $S$.
By Lemma~\ref{lem:Preliminary},
any $(p,q)$-periodic orbit in $S$ lies on $R = \Graph \zeta$.
Let $\Omega \subset S$ be the domain enclosed by the curves
$R = \Graph \zeta$ and $\hat{R} = \Graph \hat{\zeta}$.
Let $B = (\tilde{\Phi}\circ\Phi)(\Omega)$.
Let $K_2$ be the supremum of $|\det[\rmD \tilde{\Phi}]|$
in the compact set $\Tset \times [-b'_*,b'_*]$.
Let $K= 4 K_0 K_1 K_2 L (b_*')^3$.
Then, following the arguments contained in
Subsection~\ref{ssec:DifferenceOfPeriodicActions}
about the difference of periodic actions, we get that
\begin{eqnarray}
\nonumber
\Delta^{(p,q)}
& \le &
\Area[B] \le
K_1 K_2 \Area[\Omega] =
K_1 K_2 \int_\Tset \left|\hat{\zeta}(\xi) - \zeta(\xi)\right|\rmd \xi \\
\label{eq:DeltaBoundProof}
& \le &
K_1 K_2 q b^4 K_0 \rme^{-2\pi \alpha  q /|p|} \le
K \rme^{-2\pi\alpha q /|p|},
\end{eqnarray}
for all relatively prime integers $p$ and $q$ with
$1 \le |p| \le L$ and $q \ge q_*$.
We have used expression~(\ref{eq:NeishtadtForm}), $b = c^2 |p|/q \le b'_*$,
$c^2 < 4$, the bounds on the nonintegrable terms $g_j(\xi,\eta)$,
and the bounds on the Jacobians of the changes of variables
$\Phi$ and $\tilde{\Phi}$.

This ends the proof of Theorem~\ref{thm:Asymptotic} when $(m,n) = (0,1)$.

\subsection{Proof of Theorem~\ref{thm:Asymptotic}: General case}

We reduce the general case to the previous one.
We split the argument in four steps.

\noindent
\emph{Step 1: About the rational rotation numbers.}
If $C$ is a $(m,n)$-resonant RIC and $(s,r)$ is a $(p,q)$-periodic point of $g$,
then $C$ is a $(m,1)$-resonant RIC and
$(s,r)$ is a $(p',q')$-periodic point of the power map $f = g^n$,
where
\[
p' = \frac{np}{\gcd(n,q)},\qquad
q' = \frac{q}{\gcd(n,q)}.
\]
By taking the suitable lift $F$ of $f$, we can assume that
$C$ is a $(0,1)$-resonant RIC and
$(s,r)$ is a $(p'',q'')$-periodic point of $f$,
with $p''/q'' = p'/q' - m$.
That is,
\begin{equation}\label{eq:NewRotationNumber}
p'' = p' - mq' = \frac{np-mq}{\gcd(n,q)},\qquad
q'' = q' = \frac{q}{\gcd(n,q)}.
\end{equation}
If $p$ and $q$ are relatively prime integers such that $1 \le |np-mq| \le L$
and $q \ge q_*$,
$p''$ and $q''$ are relatively prime integers such that
$|p''| \le L/\gcd(n,q) \le L$ and $q'' \ge q_*/\gcd(n,q) \ge q_*/n$.

\noindent
\emph{Step 2: About the Lagrangians.}
Let $G$ and $F=G^n$ be the lifts of $g$ and $f=g^n$ we are dealing with.
If $G^* \lambda - \lambda = \rmd h$, then
\[
F^* \lambda - \lambda =
\sum_{j=0}^{n-1} \left[ \left(G^{j+1}\right)^* \lambda -
                        \left(G^j\right)^* \lambda \right] =
\sum_{j=0}^{n-1} \rmd (h \circ G^j) =
\rmd \left( \sum_{j=0}^{n-1} h \circ G^j \right),
\]
so $\ell(s_0,s_n) := h(s_0,s_1) + h(s_1,s_2) + \cdots + h(s_{n-1},s_n)$
is a Lagrangian of $f$.
This Lagrangian is well defined in a neighboorhood of the resonant RIC $C$,
because $f$ is twist on $C$.

\noindent
\emph{Step 3: About the periodic actions.}
Let $O$ be the $(p,q)$-periodic orbit of $g$ through the point $(s,r)$,
being $W^{(p,q)}[O]$ its $(p,q)$-periodic action.
Let $O''$ be the $(p'',q'')$-periodic orbit of $f$ through the same point,
being $W''^{(p'',q'')}[O'']$ its $(p'',q'')$-periodic action.
We deduce from the previous steps and a straightforward computation that
\begin{equation}\label{eq:NewPeriodicAction}
W''^{(p'',q'')}[O''] = \frac{n}{\gcd(n,q)} W^{(p,q)}[O].
\end{equation}

\noindent
\emph{Step 4: Final bound.}
The result follows directly from the bound~(\ref{eq:DeltaBoundProof})
taking into account relations~(\ref{eq:NewRotationNumber})
and~(\ref{eq:NewPeriodicAction}).
We just note that
\[
\rme^{-2 \pi \alpha q''/|p''|} =
\exp\left(-\frac{2\pi \alpha q}{|np-mq|}\right).
\]

This ends the proof of Theorem~\ref{thm:Asymptotic}.

\ack
The authors were supported in part by CUR-DIUE Grant 2014SGR504 (Catalonia)
and MINECO-FEDER Grant MTM2012-31714 (Spain).
Some of the results of this paper were obtained while A T was
in the University of Warwick as a visiting researcher.
Thanks to Vassily Gelfreich, Robert MacKay, Laurent Niederman,
S\^onia Pinto-de-Carvalho, and Tere Seara for useful and stimulating
conversations.


\section*{References}


\begin{thebibliography}{99}

\bibitem{AnderssonMelrose1977}
K.~G. Andersson and R.~B. Melrose.
\newblock The propagation of singularities along gliding rays.
\newblock {\em Invent. Math.}, 41(3):197--232, 1977.

\bibitem{ArnoldAvez1968}
V.~I. Arnol'd and A.~Avez.
\newblock {\em Ergodic Problems of Classical Mechanics}.
\newblock Translated from the French by A. Avez. W. A. Benjamin, Inc., New
  York-Amsterdam, 1968.

\bibitem{BaryshnikovZharnitsky2006}
Y.~Baryshnikov and V.~Zharnitsky.
\newblock Sub-{R}iemannian geometry and periodic orbits in classical billiards.
\newblock {\em Math. Res. Lett.}, 13(4):587--598, 2006.

\bibitem{Birkhoff1966}
G.~D. Birkhoff.
\newblock {\em Dynamical Systems}.
\newblock With an addendum by Jurgen Moser. American Mathematical Society
  Colloquium Publications, Vol. IX. American Mathematical Society, Providence,
  R.I., 1966.

\bibitem{Boyland1996}
P.~Boyland.
\newblock Dual billiards, twist maps and impact oscillators.
\newblock {\em Nonlinearity}, 9(6):1411, 1996.

\bibitem{Colin1984}
Y.~{Colin de Verdi{\`e}re}.
\newblock Sur les longueurs des trajectoires p{\'e}riodiques d'un billard.
\newblock In {\em South {R}hone seminar on geometry, {III} ({L}yon, 1983)},
  Travaux en Cours, pages 122--139. Hermann, Paris, 1984.

\bibitem{Coutinho}
L.~Coutinho.
\newblock {\em Bilhares em superficies de curvatura constante}.
\newblock PhD thesis, Universidade Federal de Minas Gerais, 2014.

\bibitem{Cyr2012}
V.~Cyr.
\newblock A number theoretic question arising in the geometry of plane curves
  and in billiard dynamics.
\newblock {\em Proc. Amer. Math. Soc.}, 140(9):3035--3040, 2012.

\bibitem{Day1947}
M.~Day.
\newblock Polygons circumscribed about closed convex curves.
\newblock {\em Trans. Amer. Math. Soc.}, 62:315--319, 1947.

\bibitem{DelshamsLlave2000}
A.~Delshams and R.~de~la Llave.
\newblock K{AM} theory and a partial justification of {G}reene's criterion for
  nontwist maps.
\newblock {\em SIAM J. Math. Anal.}, 31(6):1235--1269, 2000.

\bibitem{Douady1982}
R.~Douady.
\newblock Applications du th{\'e}or{\`e}me des tores invariantes.
\newblock Th{\`e}se, Universit{\'e} Paris VII, 1982.

\bibitem{FontichSimo1990}
E.~Fontich and C.~Sim{\'o}.
\newblock The splitting of separatrices for analytic diffeomorphisms.
\newblock {\em Ergodic Theory Dynam. Systems}, 10(2):295--318, 1990.

\bibitem{Greene1979}
J.~M. Greene.
\newblock A method for determining a stochastic transition.
\newblock {\em J. Math. Phys.}, 20:1183--1201, 1979.

\bibitem{Gutkin2012}
E.~Gutkin.
\newblock Capillary floating and the billiard ball problem.
\newblock {\em Journal of Mathematical Fluid Mechanics}, 14(2):363--382, 2012.

\bibitem{GutkinKatok1995}
E.~Gutkin and A.~Katok.
\newblock Caustics for inner and outer billiards.
\newblock {\em Comm. Math. Phys.}, 173(1):101--133, 1995.

\bibitem{Innami1988}
N.~Innami.
\newblock Convex curves whose points are vertices of billiard triangles.
\newblock {\em Kodai Mathematical Journal}, 11(1):17--24, 1988.

\bibitem{Kac1966}
M.~Kac.
\newblock Can one hear the shape of a drum?
\newblock {\em Amer. Math. Monthly}, 73(4, part II):1--23, 1966.

\bibitem{KatokHasselblatt1995}
A.~Katok and B.~Hasselblatt.
\newblock {\em Introduction to the Modern Theory of Dynamical Systems},
  volume~54 of {\em Encyclopedia of Mathematics and its Applications}.
\newblock Cambridge University Press, Cambridge, 1995.
\newblock With a supplementary chapter by Katok and Leonardo Mendoza.

\bibitem{Knill1998}
O.~Knill.
\newblock On nonconvex caustics of convex billiards.
\newblock {\em Elem. Math.}, 53(3):89--106, 1998.

\bibitem{KozlovTreschev1991}
V.V. Kozlov and D.V. Treschev.
\newblock {\em Billiards: A Genetic Introduction to the Dynamics of Systems
  with Impacts}, volume~89 of {\em Translations of Mathematical Monographs}.
\newblock Amer. Math. Soc., Providence, RI, 1991.

\bibitem{Lazutkin1973}
V.~F. Lazutkin.
\newblock Existence of caustics for the billiard problem in a convex domain.
\newblock {\em Izv. Akad. Nauk SSSR Ser. Mat.}, 37:186--216, 1973.

\bibitem{MacKay1992}
R.~S. MacKay.
\newblock Greene's residue criterion.
\newblock {\em Nonlinearity}, 5(1):161--187, 1992.

\bibitem{MackayMeissPercival1984}
R.~S. MacKay, J.~D. Meiss, and I.~C. Percival.
\newblock Transport in {H}amiltonian systems.
\newblock {\em Phys. D}, 13(1-2):55--81, 1984.

\bibitem{MarviziMelrose1982}
S.~Marvizi and R.~Melrose.
\newblock Spectral invariants of convex planar regions.
\newblock {\em J. Differential Geom.}, 17(3):475--502, 1982.

\bibitem{Mather1986}
J.~N. Mather.
\newblock A criterion for the nonexistence of invariant circles.
\newblock {\em Inst. Hautes {\'E}tudes Sci. Publ. Math.}, (63):153--204, 1986.

\bibitem{Meiss1992}
J.~D. Meiss.
\newblock Symplectic maps, variational principles, and transport.
\newblock {\em Rev. Modern Phys.}, 64(3):795--848, 1992.

\bibitem{Moser1978}
J.~Moser.
\newblock Is the solar system stable?
\newblock {\em Math. Intelligencer}, 1(2):65--71, 1978/79.

\bibitem{Neishtadt1981}
A.~I. Ne{\u\i}shtadt.
\newblock Estimates in the {K}olmogorov theorem on conservation of
  conditionally periodic motions.
\newblock {\em J. Appl. Math. Mech.}, 45(6):1016--1025, 1981.

\bibitem{Neishtadt1984}
A.~I. Ne{\u\i}shtadt.
\newblock The separation of motions in systems with rapidly rotating phase.
\newblock {\em Prikl. Mat. Mekh.}, 48(2):197--204, 1984.

\bibitem{PintodeCarvalhoRamirezRos2013}
S.~{Pinto-de-Carvalho} and R.~Ram{\'i}rez-Ros.
\newblock Non-persistence of resonant caustics in perturbed elliptic billiards.
\newblock {\em Ergodic Theory Dynam. Systems}, 33(6):1876--1890, 2013.

\bibitem{Poschel1982}
J.~P{\"o}schel.
\newblock The concept of integrability on {C}antor sets for {H}amiltonian
  systems.
\newblock {\em Celestial Mech.}, 28(1-2):133--139, 1982.

\bibitem{RamirezRos2006}
R.~Ram{\'i}rez-Ros.
\newblock Break-up of resonant invariant curves in billiards and dual billiards
  associated to perturbed circular tables.
\newblock {\em Phys. D}, 214(1):78--87, 2006.

\bibitem{SiegelMoser1995}
C.~L. Siegel and J.~K. Moser.
\newblock {\em Lectures on Celestial Mechanics}.
\newblock Classics in Mathematics. Springer-Verlag, Berlin, 1995.
\newblock Translated from the German by C. I. Kalme, Reprint of the 1971
  translation.

\bibitem{Sorrentino2014}
A.~Sorrentino.
\newblock Computing {M}ather's $\beta$-function for {B}irkhoff billiards.
\newblock To appear in {\em Discrete Contin. Dynam. Systems}, 2014.

\bibitem{Tabachnikov1995_Billiards}
S.~Tabachnikov.
\newblock Billiards.
\newblock {\em Panor. Synth.}, (1):vi+142, 1995.

\bibitem{Tabachnikov1995}
S.~Tabachnikov.
\newblock On the dual billiard problem.
\newblock {\em Adv. Math.}, 115(2):221--249, 1995.

\bibitem{Tabachnikov2002}
S.~Tabachnikov.
\newblock Dual billiards in the hyperbolic plane.
\newblock {\em Nonlinearity}, 15(4):1051--1072, 2002.

\end{thebibliography}
\end{document}